# ANOVA FOR DIFFUSIONS AND ITÔ PROCESSES[1]

By Per Aslak Mykland and Lan Zhang

*University of Chicago, and Carnegie Mellon University
and University of Illinois at Chicago*

Itô processes are the most common form of continuous semimartingales, and include diffusion processes. This paper is concerned with the nonparametric regression relationship between two such Itô processes. We are interested in the quadratic variation (integrated volatility) of the residual in this regression, over a unit of time (such as a day). A main conceptual finding is that this quadratic variation can be estimated almost as if the residual process were observed, the difference being that there is also a bias which is of the same asymptotic order as the mixed normal error term.

The proposed methodology, "ANOVA for diffusions and Itô processes," can be used to measure the statistical quality of a parametric model and, nonparametrically, the appropriateness of a one-regressor model in general. On the other hand, it also helps quantify and characterize the trading (hedging) error in the case of financial applications.

**1. Introduction.** We consider the regression relationship between two stochastic processes $\Xi_t$ and $S_t$,

$$(1.1) \qquad d\Xi_t = \rho_t \, dS_t + dZ_t, \qquad 0 \le t \le T,$$

where $Z_t$ is a residual process. We suppose that the processes $S_t$ and $\Xi_t$ are observed at discrete sampling points $0 = t_0 < \cdots < t_k = T$. With the advent of high frequency financial data, this type of regression has been a topic of growing interest in the literature; see Section 2.4. The processes $\Xi_t$ and $S_t$ will be Itô processes, which are the most commonly used type of continuous

Received March 2002; revised September 2005.

[1]Supported in part by NSF Grants DMS-99-71738 (Mykland) and DMS-02-04639 (Mykland and Zhang).

*AMS 2000 subject classifications.* Primary 60G44, 62M09, 62M10, 91B28; secondary 60G42, 62G20, 62P20, 91B84.

*Key words and phrases.* ANOVA, continuous semimartingale, statistical uncertainty, goodness of fit, discrete sampling, parametric and nonparametric estimation, small interval asymptotics, stable convergence, option hedging.

semimartingale. Diffusions are a special case of Itô processes. Definitions are made precise in Section 2.1 below. The differential in (1.1) is that of an Itô stochastic integral, as defined in Chapter I.4.d of [33] or Chapter 3.2 of [34].

Our purpose is to assess nonparametrically what is the smallest possible residual sum of squares in this regression. Specifically, for two processes $X_t$ and $Y_t$, denote the *quadratic covariation* between $X_t$ and $Y_t$ on the interval $[0, T]$ by

$$(1.2) \qquad \langle X, Y \rangle_T = \lim_{\max t_{i+1} - t_i \downarrow 0} \sum_i (X_{t_{i+1}} - X_{t_i})(Y_{t_{i+1}} - Y_{t_i}),$$

where $0 = t_0 < \cdots < t_k = T$. (This object exists by Definition I.4.45 or Theorem I.4.47 in [33], pages 51–52, and similar statements in [34] and [38].) In particular, $\langle Z, Z \rangle_T$—the *quadratic variation* of $Z_t$—is the sum of squares of the increments of the process $Z$ under the idealized condition of continuous observation. We wish to estimate, from discrete-time data,

$$(1.3) \qquad \min_\rho \langle Z, Z \rangle_T,$$

where the minimum is over all adapted regression processes $\rho$.

An important motivating application for the system (1.1) is that of statistical risk management in financial markets. We suppose that $S_t$ and $\Xi_t$ are the discounted values of two securities. At each time $t$, a financial institution is short one unit of the security represented by $\Xi$, and at the same time seeks to offset as much risk as possible by holding $\rho_t$ units of security $S$. $Z_t$, as given by (1.1), is then the gain/loss up to time $t$ from following this "risk-neutral" procedure. In a complete (idealized) financial market, $\min_\rho \langle Z, Z \rangle$ is zero; in an incomplete market, $\min_\rho \langle Z, Z \rangle$ quantifies the unhedgeable part of the variation in asset $\Xi$, when one adopts the best possible strategy using only asset $S$. And this lower bound (1.3) is the target that a risk management group wants to monitor and control.

The statistical importance of (1.3) is this. Once you know how to estimate (1.3), you know how to assess the goodness of fit of any given estimation method for $\rho_t$. You also know more about the appropriateness of a one-regressor model of the form (1.1). We return to the goodness of fit questions in Section 4. A model example of both statistical and financial importance is given in Section 2.2.

To discuss the problem of estimating (1.3), consider first how one would find this quantity if the processes $\Xi$ and $S$ were continuously observed. Note that from (1.1), one can write

$$(1.4) \qquad \langle Z, Z \rangle_t = \langle \Xi, \Xi \rangle_t + \int_0^t \rho_u^2 \, d\langle S, S \rangle_u - 2 \int_0^t \rho_u \, d\langle \Xi, S \rangle_u$$

([33], I.4.54, page 55). Here $d\langle X, Y \rangle_t$ is the differential of the process $\langle X, Y \rangle_t$ with respect to time. We shall typically assume that $\langle X, Y \rangle_t$ is absolutely continuous as a function of time (for any realization). It is easy to see that
the solution in $\rho_t$ to the problem (1.3) is uniquely given by

$$\tag{1.5} \rho_t = \frac{d\langle \Xi, S \rangle_t}{d\langle S, S \rangle_t} = \frac{\langle \Xi, S \rangle_t'}{\langle S, S \rangle_t'},$$

where $\langle \Xi, S \rangle_t'$ is the derivative of $\langle \Xi, S \rangle_t$ with respect to time. Apart from its
statistical significance, in financial terms $\rho$ is the hedging strategy associated
with the minimal martingale measure (see, e.g., [16] and [41]).

The problem (1.3) then connects to an ANOVA, as follows. Let $Z_t$ be
the residual in (1.1) for the optimal $\rho_t$, so that the quantity in (1.3) can be
written simply as $\langle Z, Z \rangle_T$. In analogy with regular regression, substituting
(1.5) into (1.4) gives rise to an ANOVA decomposition of the form

$$\tag{1.6} \underbrace{\langle \Xi, \Xi \rangle_T}_{\text{total SS}} = \underbrace{\int_0^T \rho_u^2 \, d\langle S, S \rangle_u}_{\text{SS explained}} + \underbrace{\langle Z, Z \rangle_T}_{\text{RSS}},$$

where "SS" is the abbreviation for (continuous) "sum of squares," and "RSS"
stands for "residual sum of squares." Under continuous observation, there-
fore, one solves the problem (1.3) by using the $\rho$ and $Z$ defined above. Our
target of inference, $\langle Z, Z \rangle_t$, would then be observable. Discreteness of obser-
vation, however, creates the need for inference.

The main theorems in the current paper are concerned with the asymp-
totic behavior of the estimated RSS, as more discrete observations are avail-
able within a fixed time window. There will be some choice in how to select
the estimator $\widehat{\langle Z, Z \rangle}_t$. We consider a class of such estimators $\widehat{\langle Z, Z \rangle}_t$. No
matter which of our estimators is used, we get the decomposition

$$\tag{1.7} \widehat{\langle Z, Z \rangle}_t - \langle Z, Z \rangle_t \approx \text{bias}_t + ([Z, Z]_t - \langle Z, Z \rangle_t)$$

to first-order asymptotically, where $[Z, Z]$ is the sum of squares of the incre-
ments of the (unseen) process $Z$ at the sampling points, $[Z, Z]_t = \sum_i (Z_{t_{i+1}} - Z_{t_i})^2$; see the definition (2.4) below.

A primary conceptual finding in (1.7) is the clear cut effect of the two
sources behind the asymptotics. The form of the bias depends only on the
choice of the estimator for the quadratic variation. On the other hand, the
variation component is common for all the estimators under study; it comes
only from the discretization error in discrete time sampling.

It is worthwhile to point out that our problem is only tangentially re-
lated to that of estimating the regression coefficient $\rho_t$. This is in the sense
that the asymptotic behavior of nonparametric estimators of $\rho_t$ does not
directly imply anything about the behavior of estimators of $\langle Z, Z \rangle_T$. To il-
lustrate this point, note that the convergence rates are not the same for the
two types of estimators $[O_p(n^{-1/4})$ vs. $O_p(n^{-1/2})]$, and that the variance of the estimator we use for $\rho_t$ becomes the bias in one of our estimators of $\langle Z, Z \rangle_T$ [compare equation (2.9) in Section 2.4 with Remark 1 in Section 3]. For further discussion, see Section 2.4. Depending on the goal of inference, statistical estimates $\hat{\rho}_t$ of the regression coefficient can be obtained using parametric methods, or nonparametric ones that are either local in space or in time, as discussed in Sections 2.4 and 4.1 and the references cited in these sections. In addition, it is also common in financial contexts to use calibration ("implied quantities," see Chapters 11 and 17 in [28]).

The organization is as follows: in Section 2 we establish the framework for ANOVA, and we introduce a class of estimators of the residual quadratic variation. Our main results, in Section 3, provide the distributional properties of the estimation errors for RSS. See Theorems 1 and 2. In Section 4 we discuss the statistical application of the main theorems. Parametric and nonparametric estimation are compared in the context of residual analysis. The goodness of fit of a model is addressed. Broad issues, including the analysis of variation versus analysis of variance, the moderate level of aggregation versus long run, the *actual* probability distribution versus the *risk neutral* probability distribution in the derivative valuation setting, are discussed in Sections 4.4 and 4.5. After concluding in Section 5, we give proofs in Sections 6 and 7.

## 2. ANOVA for Itô processes: framework.

2.1. *Itô processes, quadratic variation and diffusions.* The assumptions and definitions in the following two subsections are used throughout the paper, sometimes without further reference. First of all, we shall be working with a fixed filtered probability space.

System Assumption I. We suppose that there is an underlying filtered probability space $(\Omega, \mathcal{F}, \mathcal{F}_t, P)_{0 \leq t \leq T}$ satisfying the "usual conditions" (see, e.g., Definition 1.3, page 2, of [33], also in [34]).

We shall then be working with Itô processes adapted to this system, as follows. Note that Markov diffusions are a special case.

Definition 1 (Itô processes). By saying that $X$ is an *Itô process*, we mean that $X$ is continuous (a.s.), $(\mathcal{F}_t)$-adapted, and that it can be represented as a smooth process plus a local martingale:

$$(2.1) \qquad X_t = X_0 + \int_0^t \widetilde{X}_u \, du + \int_0^t \sigma_u^X \, dW_u^X,$$

where $W$ is a standard $((\mathcal{F}_t), P)$-Brownian motion, $X_0$ is $\mathcal{F}_0$ measurable, and the coefficients $\widetilde{X}_t$ and $\sigma_t^X$ are predictable, with $\int_0^T |\widetilde{X}_u| \, du < +\infty$ and

$\int_0^T (\sigma_u^X)^2 \, du < +\infty$. We also write

$$(2.2) \qquad\qquad X_t = X_0 + X_t^{\mathrm{DR}} + X_t^{\mathrm{MG}}$$

as shorthand for the drift and local martingale parts of Doob–Meyer decomposition in (2.1).

A more abstract way of putting this definition is that $X_t$ is an Itô process if it is a continuous semimartingale ([33], Definition I.4.21, page 42) whose finite variation and local martingale parts, given by (2.2), satisfy that both $X_t^{\mathrm{DR}}$ and the quadratic variation $\langle X^{\mathrm{MG}}, X^{\mathrm{MG}} \rangle_t$ are absolutely continuous. Obviously, an Itô process is a special semimartingale, also in the sense of the same definition of [33].

Diffusions are normally taken to be a special case of Itô processes, where one can write $(\sigma_t^X)^2 = a(X_t, t)$ and $\widetilde{X}_t = b(X_t, t)$, and similarly in the multidimensional setting. For a description of the link, we refer to Chapter 5.1 of [34].

The Itô process definition extends to a two- or multi-dimensional process, say, $(X_t, Y_t)$, by requiring each component $X_t$ and $Y_t$ individually to be an Itô process. Obviously, $W^X$ is typically different for different Itô processes $X$. For two processes $X$ and $Y$, the relationship between $W^X$ and $W^Y$ can be arbitrary.

The quadratic variation $\langle X, X \rangle_t$ [formula (1.2)] can now be expressed in terms of the representation (2.1) by ([33], I.4.54, page 55)

$$\langle X, X \rangle_t = \int_0^t (\sigma_u^X)^2 \, du.$$

We denote by $\langle X, X \rangle_t'$ the derivative of $\langle X, X \rangle_t$ with respect to time $t$. Then $\langle X, X \rangle_t' = (\sigma_t^X)^2$, and this quantity (or its square root) is often known as *volatility* in the finance literature.

Both quadratic variation and covariation are absolutely continuous. This follows from the Itô process assumption and from the Kunita–Watanabe inequality (see, e.g., page 69 of [38]).

DEFINITION 2 (Volatility as an Itô process). Denote by $\langle X, Y \rangle_t'$ the derivative of $\langle X, Y \rangle_t$ with respect to time. We shall often suppose that $\langle X, Y \rangle_t'$ is itself an Itô process. For ease of notation, we then write its Doob–Meyer decomposition as

$$d\langle X, Y \rangle_t' = dD_t^{XY} + dR_t^{XY} = \widetilde{D}_t^{XY} \, dt + dR_t^{XY}.$$

Note that the quadratic variation of $\langle X, Y \rangle'$ is the same as $\langle R^{XY}, R^{XY} \rangle$, and that, in our notation above, $D^{XY} = (\langle X, Y \rangle')^{\mathrm{DR}}$ and $R^{XY} = (\langle X, Y \rangle')^{\mathrm{MG}}$.

2.2. *Model example*: *Adequacy of the one factor interest rate model.* A common model for the risk free short term interest rate is given by the diffusion

$$(2.3) \qquad dr_t = \mu(r_t)\,dt + \gamma(r_t)\,dW_t,$$

where $W_t$ is a Brownian motion. For example, in [43], one takes $\mu(r) = \kappa(\alpha - r)$, and $\gamma(r) = \gamma =$ constant, while Cox, Ingersoll and Ross [6] uses the same function $\mu$, but $\gamma(r) = \gamma r^{1/2}$. For more such models, and a brief financial introduction, see, for example, [28]. A discussion and review of estimation methods is given by Fan [14].

One of the implications of this so-called one factor model is the following. Suppose $S_t$ and $\Xi_t$ are the values of two zero coupon government bonds with different maturities. Financial theory then predicts that, until maturity, $S_t = f(r_t, t)$ and $\Xi_t = g(r_t, t)$ for two functions $f$ and $g$ (see [28], Chapter 21, for details and functional forms). Under this model, therefore, the relationship $d\Xi_t = \rho_t\,dS_t$ holds from time zero until the maturity of the shorter term bond. It is easy to see that $d\langle S, S\rangle_t = f'_r(r_t, t)^2\gamma(r_t)^2\,dt$ and $d\langle \Xi, S\rangle_t = f'_r(r_t, t)g'_r(r_t, t)\gamma(r_t)^2\,dt$, with $\rho_t = d\langle \Xi, S\rangle_t/d\langle S, S\rangle_t = g'_r(r_t, t)/f'_r(r_t, t)$. Here $f'_r$ is the derivative of $f$ with respect to $r$, and similarly for $g'_r$.

The one factor model is only an approximation, and to assess the adequacy of the model, one would now wish to estimate $\langle Z, Z\rangle_T$ over different time intervals. This provides insight into whether it is worthwhile to use a one-factor model at all. If the conclusion is satisfactory, one can estimate the quantites $\mu$ and $\gamma$ (and, hence, $f$ and $g$) with parametric or nonparametric methods (see Sections 2.4 and 4.1), and again, use our methods to assess the fit of the specific model, for example, as discussed in Section 4.1.

2.3. *Finitely many data points.* We now suppose that we observe processes, in particular, $S_t$ and $\Xi_t$, on a finite set (partition, grid) $\mathcal{G} = \{0 = t_0 < t_1 < \cdots < t_k = T\}$ of time points in the interval $[0, T]$. We take the time points to be nonrandom, but possibly irregularly spaced. Note that this also covers the case where the $t_i$ are random but independent of the processes we seek to observe, so that one can condition on $\mathcal{G}$ to get back to irregular but nonrandom spacing.

DEFINITION 3 (Observed quadratic variation). For two Itô processes $X$ and $Y$ observed on a grid $\mathcal{G}$,

$$(2.4) \qquad [X, Y]_t^{\mathcal{G}} = \sum_{t_{i+1} \le t} (\Delta X_{t_i})(\Delta Y_{t_i}),$$

where $\Delta X_{t_i} = X_{t_{i+1}} - X_{t_i}$. When there is no ambiguity, we use $[X, Y]_t$ for $[X, Y]_t^{\mathcal{G}}$, or $[X, Y]_t^{(n)}$ in case of a sequence $\mathcal{G}_n$.

Note that this is not the same as the usual definition of $[X, Y]_t$ for Itô processes. *We use $\langle X, Y \rangle_t$ to refer to a continuous process, while $[X, Y]_t$ refers to a (càdlàg) process which only changes values at the partition points $t_i$.*

Since the results of this paper rely on asymptotics, we shall take limits with the help of a sequence of partitions $\mathcal{G}_n = \{0 = t_0^{(n)} < t_1^{(n)} < \cdots < t_{k_n}^{(n)} = T\}$. As $n \to +\infty$, we let $\mathcal{G}_n$ become dense in $[0, T]$, in the sense that the mesh

$$(2.5) \qquad \delta^{(n)} = \max_t |\Delta t_i^{(n)}| \to 0.$$

Here $\Delta t_i^{(n)} = t_{i+1}^{(n)} - t_i^{(n)}$. In other words, the mesh is the maximum distance between the $t_i^{(n)}$'s. On the other hand, $T$ remains fixed (except briefly in Section 4.5).

In this case, $[X, Y] = [X, Y]^{\mathcal{G}_n}$ converges to $\langle X, Y \rangle$ uniformly in probability; see [33], Theorem I.4.47, page 52, and [38], Theorem II.23, page 68. More is true; see Section 5 of [32], and (our) Sections 2.8 and 6 below.

Note that, under (2.5), $k_n = |\mathcal{G}_n| \to \infty$. It is often convenient to consider the average distance between successive observation points,

$$(2.6) \qquad \overline{\Delta t}^{(n)} = \frac{T}{k_n};$$

see Assumption A(i) below.

2.4. *The regression problem, and the estimation of $\rho_t$.* The processes in (1.1) will be taken to satisfy the following.

SYSTEM ASSUMPTION II. We let $\Xi$ and $S$ be Itô processes. We assume that

$$(2.7) \qquad \inf_{t \in [0, T]} \langle S, S \rangle_t' > 0 \qquad \text{almost surely.}$$

This assumption assures that $\rho_t$, given by (1.5), is well defined under (2.7) by the Kunita–Watanabe inequality.

As noted in the Introduction, under continuous observation of $\Xi_t$ and $S_t$, one can also directly observe the optimal $\rho_t$ and $Z_t$. Our target of inference, $\langle Z, Z \rangle_t$, would then be observable. Discreteness of observation, however, creates the need for inference.

In a noncontinuous world, where $\Xi$ and $S$ can only be observed over grid times, the most straightforward estimator of $\rho$ is

$$(2.8) \qquad \hat{\rho}_t = \frac{\widehat{\langle \Xi, S \rangle_t'}}{\widehat{\langle S, S \rangle_t'}} = \frac{[\Xi, S]_t - [\Xi, S]_{t - h_n}}{[S, S]_t - [S, S]_{t - h_n}}.$$

For simplicity, this estimator is the one we shall use in the following. The results easily generalize when more general kernels are used. Note that we have to use a smoothing bandwidth $h_n$. There will naturally be a tradeoff between $h_n$ and $\overline{\Delta t}^{(n)}$. As we now argue, this typically results in $h_n = O((\overline{\Delta t}^{(n)})^{1/2})$.

Asymptotics for estimators of the form $\widehat{\langle X, Y \rangle}'_t = ([X, Y]_t - [X, Y]_{t-h_n})/h_n$ and, hence, for $\hat{\rho}_t$, are given by Foster and Nelson [17] and Zhang [44]. Let $\overline{\Delta t}^{(n)}$ be the average observation interval, assumed to converge to zero. If $\langle X, Y \rangle'_t$ is an Itô process with nonvanishing volatility, then it is optimal to take $h_n = O((\overline{\Delta t}^{(n)})^{1/2})$, and $(\overline{\Delta t}^{(n)})^{1/4}(\widehat{\langle X, Y \rangle}'_t - \langle X, Y \rangle'_t)$ converges in law (for each fixed $t$) to a (conditional on the data) normal distribution with mean zero and random variance. (The mode of convergence is the same as in Proposition 1.) The asymptotic distributions are (conditionally) independent for different times $t$. If $\langle X, Y \rangle'_t$ is smooth, on the other hand, the rate becomes $(\overline{\Delta t}^{(n)})^{1/3}$ rather than $(\overline{\Delta t}^{(n)})^{1/4}$, and the asymptotic distribution contains both bias and variance.

The same applies to the estimator $\hat{\rho}_t$. In the case when $S$ and $\Xi$ have a diffusion component, the estimator has (random) asymptotic variance

$$(2.9) \qquad V_{\hat{\rho}-\rho}(t) = \frac{\langle \rho, \rho \rangle'_t}{3c} + cH'(t)\left( \frac{\langle \Xi, \Xi \rangle'_t}{\langle S, S \rangle'_t} - \rho_t^2 \right)$$

whenever $h_n/(\overline{\Delta t}^{(n)})^{1/2} \to c \in (0, \infty)$; see [44]. $H(t)$ is defined in Section 2.6 below.

The scheme given in (2.8) is only one of many for estimating $\langle X, Y \rangle'_t$ by using methods that are local in time. In particular, Genon-Catalot, Laredo and Picard [21] use wavelets for this purpose and determine rates of convergence and limit distributions under the assumption that $\langle X, Y \rangle'_t$ is deterministic and has smoothness properties.

Other important literature in this area seeks to estimate $\langle X, Y \rangle'_t$ as a function of the underlying state variables by methods that are local in space; see, in particular, [15, 26, 31]. The typical setup is that $U = (X, Y, \ldots)$ is a Markov process, so that $\langle X, Y \rangle'_t = f(U_t)$ for some function $f$, and the problem is to estimate $f$. If all coefficients in the Markov diffusion are smooth of order $s$, and subject to regularity conditions, the function $f$ can be estimated with a rate of convergence of $(\overline{\Delta t}^{(n)})^{s/(1+2s)}$.

The convergence obtained for the estimator of $f$ under Markov assumptions is considerably faster than what can be obtained for (2.8). It does, however, rely on stronger (Markov) assumptions than the ones (Itô processes) that we shall be working with in this paper. Since we shall only be interested in $\rho_t$ as a (random) function of time, our development does not require a Markov specification and, in particular, does not require full knowledge of what potential state variables might be.

Of course, this is just a subset of the literature for estimation of Markov diffusions. See Section 4.1 for further references.

We emphasize that the general ANOVA approach in this paper can be carried out with other schemes for estimating $\rho_t$ than the one given in (2.8). We have seen this as a matter of fine tuning and, hence, beyond the scope of this paper. This is because Theorems 1 and 2 achieve the same rate of convergence as the one obtained in Proposition 1.

2.5. *Estimation schemes for the residual quadratic variation $\langle Z, Z \rangle_t$.* We now return to the estimation of the quadratic variation $\langle Z, Z \rangle$ of residuals. Given the discrete data of $(\Xi, S)$, there are different methods to estimate the residual variation.

One scheme is to start with model (1.1). For a fixed grid $\mathcal{G}$, one first estimates $\Delta Z_{t_i}$ through the relation $\Delta \widehat{Z}_{t_i} = \Delta \Xi_{t_i} - \hat{\rho}_{t_i}(\Delta S_{t_i})$, where all increments are from time $t_i$ to $t_{i+1}$, and then obtains the quadratic variation (q.v. hereafter) of $\widehat{Z}$. This gives an estimator of $\langle Z, Z \rangle$ as

$$(2.10) \qquad [\widehat{Z}, \widehat{Z}]_t = \sum_{t_{i+1} \le t} \left( \Delta \widehat{Z}_{t_i} \right)^2 = \sum_{t_{i+1} \le t} \left[ \Delta \Xi_{t_i} - \hat{\rho}_{t_i}(\Delta S_{t_i}) \right]^2,$$

where the notation of square brackets (discrete time-scale q.v.) is invoked, since $\Delta \widehat{Z}_{t_i}$ is the increment over discrete times.

Alternatively, one can directly analyze the ANOVA version (1.6) of the model, where $d\langle Z, Z \rangle_t = d\langle \Xi, \Xi \rangle_t - \rho_t^2 \, d\langle S, S \rangle_t$. This yields a second estimator of $\langle Z, Z \rangle_t$,

$$(2.11) \qquad \widehat{\langle Z, Z \rangle}_t^{(1)} = \sum_{t_{i+1} \le t} \left[ (\Delta \Xi_{t_i})^2 - \hat{\rho}_{t_i}^2 (\Delta S_{t_i})^2 \right].$$

In general, any convex combination of these two,

$$(2.12) \qquad \widehat{\langle Z, Z \rangle}_t^{(\alpha)} = (1 - \alpha)[\widehat{Z}, \widehat{Z}]_t + \alpha \widehat{\langle Z, Z \rangle}_t^{(1)},$$

would seem like a reasonable way to estimate $\langle Z, Z \rangle_t$, and this is the class of estimators that we shall consider. Particular properties will be seen to attach to $\widehat{\langle Z, Z \rangle}_t^{(1/2)}$, which we shall also denote by $\widetilde{\langle Z, Z \rangle}_t$. For a start, it is easy to see that

$$(2.13) \qquad \widetilde{\langle Z, Z \rangle}_t = [\Xi, \widehat{Z}]_t.$$

Note that (2.13) also has a direct motivation from the continuous model. Since $\langle S, Z \rangle_t = 0$, (1.1) yields that $\langle \Xi, Z \rangle_t = \langle Z, Z \rangle_t$.

We establish the statistical properties of the estimator $\widehat{\langle Z, Z \rangle}_t^{(\alpha)}$ and, in particular, those of $([\widehat{Z}, \widehat{Z}]$ and $\widetilde{\langle Z, Z \rangle}_t)$ in Section 3. Asymptotic properties are naturally studied with the help of small interval asymptotics.

2.6. *Paradigm for asymptotic operations.* The asymptotic property of the estimation error is considered under the following paradigm.

ASSUMPTION A (*Quadratic variation of time*). For each $n \in N$, we have a sequence of nonrandom partitions $\mathcal{G}_n = \{t_i^{(n)}\}$, $\Delta t_i^{(n)} = t_{i+1}^{(n)} - t_i^{(n)}$. Let $\max_i (\Delta t_i^{(n)}) = \delta^{(n)}$. Suppose that:

(i) $\delta^{(n)} \to 0$ as $n \to \infty$, and $\delta^{(n)}/\overline{\Delta t}^{(n)} = O(1)$.

(ii) $H_{(n)}(t) = \frac{\sum_{t_{i+1}^{(n)} \leq t} (\Delta t_i^{(n)})^2}{\overline{\Delta t}^{(n)}} \to H(t)$ as $n \to \infty$.

(iii) $H(t)$ is continuously differentiable.

(iv) The bandwidth $h_n$ satisfies $\frac{\sqrt{\overline{\Delta t}^{(n)}}}{h_n} \to c$, where $0 < c < \infty$.

(v) $[H_{(n)}(t) - H_{(n)}(t - h_n)]/h_n \to H'(t)$, where the convergence is uniform in $t$.

When the partitions are evenly spaced, $H(t) = t$ and $H'(t) = 1$. In the more general case, the left-hand side of (ii) is bounded by $t\delta^{(n)}/\overline{\Delta t}^{(n)}$, while the left-hand side of (v) is bounded by $\delta^{(n)^2}/(\overline{\Delta t}^{(n)} h) + \delta^{(n)}/\overline{\Delta t}^{(n)}$. In all our results, $h$ is eventually bigger than $\overline{\Delta t}^{(n)}$ and, hence, both the left-hand sides are bounded because of (i). The assumptions in (ii) and (v) are, therefore, about a unique limit point, and about interchanging limits and differentiation.

Note that we are not assuming that the grids are nested. Also, as discussed in Section 2.4, how fast $h_n$ and $\overline{\Delta t}^{(n)}$, respectively, decay has a trade-off in terms of the asymptotic variance of the estimation error in $\rho$. It is optimal to take $h_n = O(\sqrt{\overline{\Delta t}^{(n)}})$, whence Assumption A(iv). From now on, we use $h$ and $h_n$ interchangeably.

2.7. *Assumptions on the process structure.* The following assumptions are frequently imposed on the relevant Itô processes.

ASSUMPTION B($X$) (*Smoothness*). $X$ is an Itô process. Also, $\langle X, X \rangle_t'$ and $\widetilde{X}_t$ are continuous almost surely.

The addition of Assumption B to an Itô process $X$, and similar smoothness assumptions in results below, is partially due to the estimation of $\rho$, which requires stronger smoothness conditions. In some instances, Assumption B is partially also a matter of convenience in a proof and can be dropped at the cost of more involved technical arguments.

2.8. *The limit for the discretization error.* The error $\widehat{\langle Z, Z \rangle}_t - \langle Z, Z \rangle_t$ can be decomposed into bias and pure discretization error $[Z, Z]_t - \langle Z, Z \rangle_t$. We here discuss the limit result for the latter, following [32]. We first need the following.

SYSTEM ASSUMPTION III (*Description of the filtration*). There is a continuous multidimensional $P$-local martingale $\mathcal{X} = (\mathcal{X}^{(1)}, \ldots, \mathcal{X}^{(p)})$, for some finite $p$, so that $\mathcal{F}_t$ is the smallest sigma-field containing $\sigma(\mathcal{X}_s, s \leq t)$ and $\mathcal{N}$, where $\mathcal{N}$ contains all the null sets in $\sigma(\mathcal{X}_s, s \leq T)$.

The final statement in the assumption assures that the "usual conditions" ([33], page 2, [34], page 10) are satisfied. The main implication, however, is on our mode of convergence, as follows.

PROPOSITION 1 (Discretization theorem). *Let $Z$ be an Itô process for which $\int_0^T (\langle Z, Z \rangle')_t^2 \, dt < \infty$ a.s. and $\int_0^T \widetilde{Z}_t^2 \, dt < \infty$ a.s. Subject to Assumptions* A(i)–(ii) *and System Assumptions* I *and* III,

$$(\overline{\Delta t}^{(n)})^{-1/2} ([Z, Z]_t^{(n)} - \langle Z, Z \rangle_t) \xrightarrow{\mathcal{L}.stable} \int_0^t \sqrt{2H'(u)} \langle Z, Z \rangle_u' \, dW_u,$$

*where $W$ is a standard Brownian motion, independent of the underlying data process $\mathcal{X}$.*

The symbol $\xrightarrow{\mathcal{L}.stable}$ denotes *stable* convergence of the process, as defined in [39] and [1]; see also [40] and Section 2 of [32].

In the case of an equidistant grid, the result coincides with the applicable part of Theorems 5.1 and 5.5 in [32], and the proof is essentially the same (see Section 6). In abstract form, results of this type appear to go back to [40]. The Jacod and Protter result was used in financial applications by Zhang [44] and Barndorff-Nielsen and Shephard [3]. The case where $Z_t$ is observed discretely and with additive error is considered in [45] and [46].

Note that the conditions on $\langle Z, Z \rangle'$ and $\widetilde{Z}$ are the same as in the equidistant case, due to the Lipschitz continuity of $H$. Some further discussion and results are contained in Section 6.

## 3. ANOVA for diffusions: main distributional results.

3.1. *Distribution of $[\widehat{Z}, \widehat{Z}]_t - \langle Z, Z \rangle_t$.* Recall that the square bracket $[Z, Z]$ and the angled bracket $\langle Z, Z \rangle$ represent the quadratic variation of $Z$ at discrete and continuous time-scale, respectively.

THEOREM 1. *Under System Assumptions* I–II [*and, in particular, equation* (1.1)], *assume that Assumption* A *holds. Suppose that* $S$, $\Xi$, $\rho$, $\langle S, S \rangle'$, $\langle \Xi, S \rangle'$, $\langle R^{SS}, R^{SS} \rangle'$, $\langle R^{\Xi S}, R^{\Xi S} \rangle'$, *and* $\langle R^{\Xi\Xi}, R^{\Xi\Xi} \rangle'$ *are Itô processes, each satisfying Assumption* B. *Let the estimator* $[\widehat{Z}, \widehat{Z}]_t$ *be defined as in* (2.10). *Then, as* $n \to \infty$,

$$
\begin{aligned}
(3.1) \qquad (\overline{\Delta t}^{(n)})^{-1/2} &([\widehat{Z}, \widehat{Z}]_t - \langle Z, Z \rangle_t) \\
&= D_t + (\overline{\Delta t}^{(n)})^{-1/2}([Z, Z]_t^{(n)} - \langle Z, Z \rangle_t) + o_p(1)
\end{aligned}
$$

*uniformly in* $t$, *where*

$$
(3.2) \qquad D_t = \frac{1}{3c} \int_0^t \langle \rho, \rho \rangle'_u \, d\langle S, S \rangle_u + c \int_0^t H'(u) \, d\langle Z, Z \rangle_u.
$$

REMARK 1. The consequence of Theorem 1 is that the quantity in (3.1) converges in law (stably) to

$$
D_t + \int_0^t \sqrt{2H'(u)} \langle Z, Z \rangle'_u \, dW_u;
$$

the $o_p(1)$ term goes away by Lemma VI.3.31, page 352 in [33].

Note that $D_t$ in (3.2) can be expressed as $D_t = \int_0^t V_{\hat\rho - \rho}(u) \, d\langle S, S \rangle_u$, where $V_{\hat\rho - \rho}(t)$ is the asymptotic variance of $\hat\rho_t - \rho_t$; see (2.9) or [44]. Hence, the (random) variance term for $\hat\rho$ becomes a bias term for $[\widehat{Z}, \widehat{Z}]$. This is intuitively natural since the $\hat\rho_t$ are asymptotically independent for different $t$.

Theorem 1, together with Proposition 1, says that the estimator $[\widehat{Z}, \widehat{Z}]_t$ converges to $\langle Z, Z \rangle_t$ at the order of the square root of the average sampling interval. In the limit the error term consists of a nonnegative bias $D_t$, due to the estimation uncertainty $[\widehat{Z}, \widehat{Z}] - [Z, Z]$, and a mixture Gaussian, due to the discretization $[Z, Z]_t - \langle Z, Z \rangle_t$. The nonnegativeness of the asymptotic bias occurs because the q.v.'s $(\langle \rho, \rho \rangle, \langle S, S \rangle, \langle Z, Z \rangle)$ are nondecreasing processes. Furthermore, (3.2) displays a bias–bias tradeoff; thus, an optimal $c$ for smoothing can be reached to minimize the asymptotic bias, though we have not investigated the effect of having a random $c$. The discretization term is independent of the smoothing factor.

### 3.2. *Distribution of* $\widetilde{\langle Z, Z \rangle}_t - \langle Z, Z \rangle_t$.

THEOREM 2. *Under System Assumptions* I–II, *assume that Assumption* A *holds. Also assume each of the following processes exists, and is an Itô-process satisfying Assumption* B: $\Xi$, $S$, $\rho$, $\langle \Xi, S \rangle'$, $\langle S, S \rangle'$, $\langle R^{SS}, R^{SS} \rangle'$, $\langle R^{\Xi S}, R^{SS}\rangle'$ and $\langle R^{\Xi S}, R^{\Xi S}\rangle'$. Also suppose that the processes $\langle \Xi, \rho \rangle'$ and $\langle S, \rho \rangle'$ are continuous. Then, uniformly in $t$,

$$
\begin{aligned}
(3.3) \quad & (\overline{\Delta t}^{(n)})^{-1/2}(\widetilde{\langle Z, Z \rangle}_t - \langle Z, Z \rangle_t) \\
& = \frac{1}{2c} \int_0^t \langle \Xi, S \rangle'_u \, d\rho_u \\
& \quad + (\overline{\Delta t}^{(n)})^{-1/2}([Z, Z]_t^{(n)} - \langle Z, Z \rangle_t) + o_p(1).
\end{aligned}
$$

Remark 1 applies similarly.

Unlike $[\widehat{Z}, \widehat{Z}]$, the asymptotic (conditional) bias associated with $\widetilde{\langle Z, Z \rangle}_t$ does not necessarily have a positive or negative sign. Moreover, we are no longer faced with a bias–bias tradeoff due to the position of $c$ in (3.3). In this case the role of smoothing in the asymptotic bias will be discussed in Section 3.3.

3.3. *General results for the* $\widehat{\langle Z, Z \rangle}_t^{(\alpha)}$ *class of estimators.* From (2.12),

$$
\widehat{\langle Z, Z \rangle}_t^{(\alpha)} = (1 - 2\alpha)[\widehat{Z}, \widehat{Z}]_t + 2\alpha \widetilde{\langle Z, Z \rangle}_t,
$$

and it follows from the assumptions of Theorems 1 and 2 that, if one sets

$$
(3.4) \quad \mathrm{bias}_t^{(\alpha)} = \frac{\alpha}{c} \int_0^t \langle \Xi, S \rangle'_u \, d\rho_u + (1 - 2\alpha)D_t,
$$

then, as in Remark 1,

$$
\begin{aligned}
(3.5) \quad & (\overline{\Delta t}^{(n)})^{-1/2}(\widehat{\langle Z, Z \rangle}_t^{(\alpha)} - \langle Z, Z \rangle_t) \\
& = \mathrm{bias}_t^{(\alpha)} + (\overline{\Delta t}^{(n)})^{-1/2}([Z, Z]_t^{(n)} - \langle Z, Z \rangle_t) + o_p(1) \\
& \xrightarrow{\mathcal{L}\text{-}stable} \mathrm{bias}_t^{(\alpha)} + \int_0^t \sqrt{2H'(u)} \langle Z, Z \rangle'_u \, dW_u.
\end{aligned}
$$

In summary, for any linear combination of the estimators in Theorems 1 and 2, $\alpha \in [0, 1]$, the convergence in (3.5) is in law as a process, and the limiting Brownian motion $W$ is independent of the entire data process. For details of stable convergence, see the discussion and references in Section 2.8 above.

The "variance" term $(\overline{\Delta t}^{(n)})^{-1/2}([Z, Z]_t - \langle Z, Z \rangle_t)$ is the same for any estimator in the linear-combination class, and they are all asymptotically perfectly correlated. The common asymptotic, conditional variance is independent of the smoothing bandwidth. It remains unclear whether the common asymptotic variance could, perhaps, be a lower bound under the nonparametric setting (see [5] for a comprehensive discussion). This needs further investigation.

TABLE 1
*The effect of constant $\rho$ on the bias components*

| Estimator | Asymptotic bias |
|---|---|
| $[\widehat{Z}, \widehat{Z}]_t$ | $c \int_0^t H'(u)\, d\langle Z, Z\rangle_u$ |
| $\widehat{\langle Z, Z\rangle}_t^{(1)}$ | $-c \int_0^t H'(u)\, d\langle Z, Z\rangle_u$ |
| $\widetilde{\langle Z, Z\rangle}_t$ | $0$ |

For the bias, on the other hand, the estimation procedure plays an important role, as the bias term varies with $\alpha$. Also the smoothing effect enters the bias terms. From Theorems 1 and 2, excessive over-smoothing (smaller $c$) or under-smoothing (bigger $c$) can explode the bias of $\widehat{\langle Z, Z\rangle}^{(\alpha)}$, for $\alpha \neq \frac{1}{2}$, thus (conditional) bias may be minimized optimally. When $\alpha = \frac{1}{2}$, it is not obvious how to deal with bias–bias tradeoff. One might theoretically be able to reduce the bias for $\widetilde{\langle Z, Z\rangle}$ [i.e., $\widehat{\langle Z, Z\rangle}^{(1/2)}$] by choosing the smallest possible bandwidth $h$. This thought should, however, be taken with caution. It is not obvious whether the magnitude of the higher-order terms in the earlier results would remain negligible if the estimation window $h$ were to decrease faster than the order $\sqrt{\overline{\Delta t}}$.

Table 1 shows that assuming constant $\rho$, $\widetilde{\langle Z, Z\rangle}_t$ will be the best choice among the three. When $\rho$ is random, none of the estimation schemes in Section 2.5 is obviously superior to the others.

3.4. *Estimating the asymptotic distribution.* For statistical inference concerning $\langle Z, Z\rangle$, one needs, in view of the above, to estimate the asymptotic (random) bias and variance. The bias term can be obtained by substitution of estimated quantities into the relevant expressions, most generally (3.4). We shall here be concerned with the variance term in (3.5). In view of the stable convergence in Proposition 1, we therefore seek to estimate $\int_0^t 2H'(u)(\langle Z, Z\rangle_u')^2\, du$.

By a modification of Barndorff-Nielsen and Shephard [3], and also Mykland [37], one can do this by considering the fourth-order variation.

DEFINITION 4 (Observed fourth-order variation). For an Itô processes $X$ observed on a grid $\mathcal{G}$,

$$(3.6) \qquad [X, X, X, X]_t = \sum_{t_{i+1} \leq t} (\Delta X_{t_i})^4.$$

PROPOSITION 2 (Estimation of variance). *Assume the regularity conditions of Theorem 1, and let $\widehat{Z}$ be defined as in that result. Also assume*

*System Assumption* III. *Then, as* $n \to \infty$,

$$(3.7) \qquad \tfrac{2}{3}(\overline{\Delta t}^{(n)})^{-1}[\widehat{Z}, \widehat{Z}, \widehat{Z}, \widehat{Z}]_t \to \int_0^t 2H'(u)(\langle Z, Z \rangle_u')^2 \, du$$

*uniformly in probability.*

This estimate of variance can be used in connection with Sections 4.2 and 4.3 below.

The proof is given in Section 6. It can be noted from there that the same statement (3.7) would hold under weaker conditions if $\widehat{Z}$ were replaced by $Z$, as follows.

REMARK 2. Assume the conditions of Proposition 1, and also that $|\widetilde{Z}|$ and $\langle Z, Z \rangle'$ are bounded a.s. Then, as $n \to \infty$,

$$(3.8) \qquad \tfrac{2}{3}(\overline{\Delta t}^{(n)})^{-1}[Z, Z, Z, Z]_t \to \int_0^t 2H'(u)(\langle Z, Z \rangle_u')^2 \, du$$

*uniformly in probability.*

This generalizes the corresponding result at the end of page 270 in [3]. The finding in [37] is exact for small samples in the context of explicit embedding, where it follows from Bartlett identities. For another use of this methodology, see, for example, the proof of Lemma 1 in [35].

**4. Goodness of fit.** The purpose of ANOVA is to assess the goodness of fit of a regression model on the form (1.1). We here illustrate the use of Theorems 1 and 2 by considering two different questions of this type. In the first section we discuss how to assess the fit of a parametric estimator for $\rho$. Afterward, we focus on the issue of how good is the one regressor model itself, independently of estimation techniques. This is already measured by the quantity $\langle Z, Z \rangle_T$, but can be further studied by considering confidence bands for $\langle Z, Z \rangle_t$ as a process, and by an analogue to the coefficient of determination. Finally, we discuss the question of the relationship between this ANOVA and the analysis of variance that is used in the standard regression setting.

4.1. *The assessment of parametric models.* In the following we suppose that a parametric model is fit to the data, and $\rho$ is estimated as a function of the parameter. Parametric estimation of discretely observed diffusions has been studied by Genon-Catalot and Jacod [18], Genon-Catalot, Jeantheau and Laredo [19, 20], Gloter [22], Gloter and Jacod [23], Barndorff-Nielsen and Shephard [2], Bibby, Jacobsen and Sørensen [4], Elerian, Siddhartha and Shephard [13], Jacobsen [29], Sørensen [42] and Hoffmann [27]. This is, of course, only a small sample of the literature available. Also, these references only concern the type of asymptotics considered in this paper, where $[0, T]$ is fixed and $\overline{\Delta t} \to 0$, and there is also a substantial literature on the case where $\overline{\Delta t}$ is fixed and $T \to \infty$.

In Section 3 we have studied the nonparametric estimate $\widehat{\langle Z, Z \rangle}_T^{(\alpha)}$ for the residual sum of squares. We here compare $\widehat{\langle Z, Z \rangle}_T^{(\alpha)}$ to its parametric counterpart to see how good the parametric model is in capturing the true regression of $\Xi$ on $S$.

Specifically, we suppose that data from the multidimensional process $X_t$ is observed at the grid points. $X_t$ has among its components at least $S_t$ and $\Xi_t$, and possibly also other processes. The parametric model is of the form $P_{\theta, \psi}$, $\theta \in \Theta$, $\psi \in \Psi$, where the modeling is such that diffusion coefficients are functions of $\theta$, while drift coefficients can be functions of both $\theta$ and $\psi$. It is thus reasonable to suppose that as $\overline{\Delta t} \to 0$, $\hat{\theta}$ converges in probability to a nonrandom parameter value $\theta_0$, and that

$$(\overline{\Delta t})^{-1/2}(\hat{\theta} - \theta_0) \to \eta N(0, 1)$$

in law stably, where $\eta$ is a function of the data and the $N(0, 1)$ term is independent of the data. (For conditions under which this occurs, consult, e.g., the references cited above.) $\theta_0$ is the true value of the parameter if the model does contain the true probability, but is otherwise also taken to be a defined parameter.

Under $P_{\theta, \psi}$, the regression coefficient $\rho_t$ is of the form $\beta_t(\theta)$. Most commonly, $\beta_t(\theta) = b(X_t; \theta)$ for a nonrandom functional $b$.

We now ask whether the true regression coefficient can be correctly estimated with the model at hand. In other words, we wish to test the null hypothesis $H_0$ that $\beta_t(\theta_0) = \rho_t$.

For the ANOVA analysis, define the theoretical residual by

$$dV_t = d\Xi_t - \beta_t(\theta_0)\, dS_t, \qquad V_0 = 0,$$

and the observed one by

$$\Delta \widehat{V}_{t_i} = \Delta \Xi_{t_i} - \beta_{t_i}(\hat{\theta})\Delta S_{t_i}, \qquad \widehat{V}_0 = 0.$$

Under the null hypothesis $\langle V, V \rangle = \langle Z, Z \rangle$, and so a natural test statistic is of the form

$$U = (\overline{\Delta t})^{-1/2}([\widehat{V}, \widehat{V}]_T - \widehat{\langle Z, Z \rangle}_T^{(\alpha)}).$$

We now derive the null distribution for $U$, using the results above.

As an intermediate step, define the discretized theoretical residual

$$\Delta V_{t_i}^d = \Delta \Xi_{t_i} - \beta_{t_i}(\theta_0)\Delta S_{t_i}, \qquad V_0^d = 0.$$

Subject to obvious regularity conditions,

$$
\begin{aligned}
[\widehat{V}, \widehat{V}]_T - [V^d, V^d]_T = &- \sum_{t_{i+1} \leq t} (\beta_{t_i}(\hat{\theta}) - \beta_{t_i}(\theta_0)) \Delta V_{t_i}^d \Delta S_{t_i} \\
&+ \sum_{t_{i+1} \leq t} (\beta_{t_i}(\hat{\theta}) - \beta_{t_i}(\theta_0))^2 (\Delta S_{t_i})^2 \\
= &-2(\hat{\theta} - \theta_0) \sum_{t_{i+1} \leq t} \frac{\partial \beta_{t_i}}{\partial \theta}(\theta_0) \Delta V_{t_i}^d \Delta S_{t_i} + O_p(\overline{\Delta t}) \\
= &-2(\hat{\theta} - \theta_0) \int_0^T \frac{\partial \beta_t}{\partial \theta}(\theta_0) \, d\langle V, S \rangle_t + O_p(\overline{\Delta t}).
\end{aligned}
$$

Also, under the conditions in Proposition 3 in Section 6 (cf. also the proof of Proposition 2 in the same section), $[V^d, V^d]_T = [V, V]_T + o_p(\overline{\Delta t}^{1/2})$ as $\overline{\Delta t} \to 0$ [since $\langle V^d, V^d \rangle_t = \langle V, V \rangle_t + o_p(\overline{\Delta t}^{1/2})$, and $\langle V^d, V^d \rangle_t' \approx \langle V^d, V \rangle_t' \approx \langle V, V \rangle_t'$].

Hence, under the conditions of Theorem 1 or 2,

$$
\begin{aligned}
U = &-2(\overline{\Delta t})^{-1/2}(\hat{\theta} - \theta_0) \int_0^T \frac{\partial \beta_t}{\partial \theta}(\theta_0) \, d\langle V, S \rangle_t \\
&+ (\overline{\Delta t})^{-1/2} \{([V, V]_T - [Z, Z]_T)\} \\
&- \operatorname{bias}_T^{(\alpha)} + o_p(1),
\end{aligned}
$$

where $\operatorname{bias}_T^{(\alpha)}$ has the same meaning as in Section 3. If the null hypothesis is satisfied, therefore,

$$
U \to N(0, 1) \times 2\eta \int_0^T \frac{\partial \beta_t}{\partial \theta}(\theta_0) \, d\langle V, S \rangle_t - \operatorname{bias}_T^{(\alpha)}
$$

in law stably. The variance and bias can be estimated from the data. This, then, provides the null distribution for $U$.

Another approach is to use $U$ to measure how close the parametric residual $\langle V, V \rangle$ is to the lower bound $\langle Z, Z \rangle$. To first order,

$$
\begin{aligned}
(\overline{\Delta t})^{1/2} U \overset{P}{\longrightarrow} &\langle V, V \rangle_T - \langle Z, Z \rangle_T \\
= &\int_0^T (\beta_t(\theta_0) - \rho_t)^2 \, d\langle S, S \rangle_t.
\end{aligned}
$$

The behavior of $U - (\overline{\Delta t})^{-1/2}(\langle V, V \rangle_T - \langle Z, Z \rangle_T)$ depends on the joint limiting distribution of $([V, V]_T - \langle V, V \rangle_T) - ([Z, Z]_T - \langle Z, Z \rangle_T)$ and $(\overline{\Delta t})^{-1/2}(\hat{\theta} - \theta_0)$. The former can be provided by Proposition 1 in Section 2.8 (or Section 5 of [32]), but further assumptions are needed to obtain the joint distribution. A study of this is beyond the scope of this paper.

4.2. *Confidence bands.* In addition to providing pointwise confidence intervals for $\widehat{\langle Z, Z \rangle}_t^{(\alpha)}$, we can also construct joint confidence bands for the estimated quadratic variation $\widehat{\langle Z, Z \rangle}^{(\alpha)}$ of residuals. This is possible because $\widehat{\langle Z, Z \rangle}^{(\alpha)}$ converges as a process by Theorems 1 and 2.

One proceeds as follows. As a process on $[0, T]$,

$$(\overline{\Delta t}^{(n)})^{-1/2}(\widehat{\langle Z, Z \rangle}_t^{(\alpha)} - \langle Z, Z \rangle_t) \xrightarrow{\mathcal{L}} \mathrm{bias}_t^{(\alpha)} + L_t.$$

Under all estimation schemes in the linear combination class, we have, by Theorems 1 and 2 and subsequent results on $\widehat{\langle Z, Z \rangle}_t^{(\alpha)}$,

$$L_t = \int_0^t \sqrt{2H'(u)} \langle Z, Z \rangle_u' \, dW_u,$$

where $W$ is a standard Brownian motion independent of the complete data filtration. Now condition on $\mathcal{F}_T$: by the stable convergence, $L_t$ is then a Gaussian process, with $\langle L, L \rangle_t$ nonrandom. Use the change-of-time construction of Dambis [7] and Dubins and Schwarz [10] to obtain $L_t = W^*_{\langle L, L \rangle_t}$, where $W^*$ is a new Brownian motion conditional on $\mathcal{F}_T$. It then follows that

$$\max_{0 \le t \le T} L_t = \max_{0 \le t \le 2\int_0^T H'(u)(\langle Z, Z \rangle_u')^2 \, du} W^*_t,$$

$$\min_{0 \le t \le T} L_t = \min_{0 \le t \le 2\int_0^T H'(u)(\langle Z, Z \rangle_u')^2 \, du} W^*_t.$$

Now write $L_n(t) = (\overline{\Delta t}^{(n)})^{-1/2}(\widehat{\langle Z, Z \rangle}_t^{(\alpha)} - \langle Z, Z \rangle_t) - \mathrm{bias}_t^{(\alpha)}$. We have

$$P(|L_n(t)| \le c, \text{ for all } t \in [0, T]) \to P(|L(t)| \le c, \text{ for all } t \in [0, T])$$

$$= P\left( \min_{0 \le t \le \tau} W^*_t \ge -c, \max_{0 \le t \le \tau} W^*_t \le c \right).$$

Choose $c = c_\tau$ such that

$$P\left( \min_{0 \le t \le \tau} W^*_t \ge -c_\tau, \max_{0 \le t \le \tau} W^*_t \le c_\tau \Big| \tau \right) = 1 - \alpha,$$

with $\tau = 2\int_0^T H'(u)(\langle Z \rangle_u')^2 \, du$. To find $c_\tau$, one can refer to Section 2.8 in [34] for the distributions of the running minimum and maximum of a Brownian motion. $\tau$ itself can be estimated by using Proposition 2. This completes our construction of a global confidence band.

4.3. *The coefficient of determination $R^2$.* In analogy with standard linear regression, one can define $R^2$ by

$$R_t^2 = 1 - \frac{\langle Z, Z \rangle_t}{\langle \Xi, \Xi \rangle_t}.$$

This quantity would have been observed if the whole paths of the processes $\Xi$ and $S$ had been available. If observations are on a grid, it is natural to use

$$\widehat{R}_t^2 = 1 - \frac{\widehat{\langle Z, Z \rangle}_t^{(\alpha)}}{[\Xi, \Xi]_t}.$$

Under the assumptions of Section 2, the distribution of $\widehat{R}_t^2$ can be found by

$$(\overline{\Delta t}^{(n)})^{-1/2}(\widehat{R}_t^2 - R_t^2)$$

$$= -(\overline{\Delta t}^{(n)})^{-1/2}\left[\frac{\widehat{\langle Z, Z \rangle}_t^{(\alpha)} - \langle Z, Z \rangle_t}{\langle \Xi, \Xi \rangle_t} - (1 - R_t^2)\frac{[\Xi, \Xi]_t - \langle \Xi, \Xi \rangle_t}{\langle \Xi, \Xi \rangle_t}\right]$$

$$\quad + o_p(1)$$

$$= -(\overline{\Delta t}^{(n)})^{-1/2}$$

$$\quad \times \frac{1}{\langle \Xi, \Xi \rangle_t}(([Z, Z]_t - \langle Z, Z \rangle_t) - (1 - R_t^2)([\Xi, \Xi]_t - \langle \Xi, \Xi \rangle_t))$$

$$\quad - \frac{\text{bias}_t^{(\alpha)}}{\langle \Xi, \Xi \rangle_t} + o_p(1),$$

where $\text{bias}_t^{(\alpha)}$ is the bias corresponding to the estimator $\widehat{\langle Z, Z \rangle}_t^{(\alpha)}$.

A straightforward generalization of Proposition 1 yields that $(\overline{\Delta t}^{(n)})^{-1/2} \times ([Z, Z]_t - \langle Z, Z \rangle_t, [\Xi, \Xi]_t - \langle \Xi, \Xi \rangle_t)_{0 \le t \le T}$ converges (stably) to a process with (bivariate) quadratic variation $\int_0^t g_u \, du$, where

$$g_t = 2H'(t)\begin{pmatrix} (\langle Z, Z \rangle_t')^2 & (\langle Z, \Xi \rangle_t')^2 \\ (\langle Z, \Xi \rangle_t')^2 & (\langle \Xi, \Xi \rangle_t')^2 \end{pmatrix},$$

and equation (1.1) yields that $\langle Z, \Xi \rangle_t' = \langle Z, Z \rangle_t'$. It follows that

$$(\overline{\Delta t}^{(n)})^{-1/2}(\widehat{R}_t^2 - R_t^2)$$

$$\xrightarrow{\mathcal{L}\text{-}stable} \frac{R_t^2}{\langle \Xi, \Xi \rangle_t}\int_0^t \sqrt{2H'(u)}\langle Z, Z \rangle_u' \, dW_u$$

$$\quad + \frac{1 - R_t^2}{\langle \Xi, \Xi \rangle_t}\int_0^t \sqrt{2H'(u)[(\langle \Xi, \Xi \rangle_u')^2 - (\langle Z, Z \rangle_u')^2]} \, dW_u^* - \frac{\text{bias}_t^{(\alpha)}}{\langle \Xi, \Xi \rangle_t},$$

where $W$ and $W^*$ are independent Brownian motions. For fixed $t$, the limit is conditionally normal, with mean $-\text{bias}_t^{(\alpha)}/\langle \Xi, \Xi \rangle_t$ and variance

$$\frac{1}{\langle \Xi, \Xi \rangle_t^2}R_t^4\int_0^t 2H'(u)(\langle Z, Z \rangle_u')^2 \, du$$

$$+ \frac{1}{\langle \Xi, \Xi \rangle_t^2} (1 - R_t^2)^2 \int_0^t 2H'(u)[(\langle \Xi, \Xi \rangle_u')^2 - (\langle Z, Z \rangle_u')^2] \, du,$$

which can be readily estimated using the device discussed in Section 3.4.

4.4. *Variance versus variation*: *Which ANOVA?* The formulation of (1.6) is in terms of quadratic variation. This raises the question of how our analysis relates to the traditional meaning of ANOVA, namely, a decomposition of variance. There are several answers to this, some concerning the broad setting provided by model (1.1), and they are discussed presently. More specific structure is provided by financial applications, and a discussion is provided in Section 4.5.

In model (1.1), the variation in $Z$ can come from both the drift and martingale components. As in (2.2),

$$(4.1) \qquad\qquad Z_t = Z_0 + Z_t^{\mathrm{DR}} + Z_t^{\mathrm{MG}}.$$

Our analysis concerns most directly the variation in $Z_t^{\mathrm{MG}}$, in that $\mathrm{var}(Z_t^{\mathrm{MG}}) = E(\langle Z, Z \rangle_t)$, where it should be noted that $\langle Z, Z \rangle_t = \langle Z^{\mathrm{MG}}, Z^{\mathrm{MG}} \rangle_t$. Hence, if the $Z^{\mathrm{DR}}$ term is identically zero, the analysis of variation is an exact analysis of variance, in terms of expectations. The quadratic variation, however, is also a more relevant measure of variation for the data that were actually collected. The Dambis [7] and Dubins and Schwarz [10] representation yields that $Z_t^{\mathrm{MG}} = V_{\langle Z, Z \rangle_t}$, where $V$ is a standard Brownian motion on a different time scale. Therefore, $\langle Z, Z \rangle_t = \langle Z^{\mathrm{MG}}, Z^{\mathrm{MG}} \rangle_t$ contains information about the actual amount of variation that has occurred in the process $Z_t^{\mathrm{MG}}$. Using the quadratic variation is, in this sense, analogous to using observed information in a likelihood setting (see, e.g., [12]). The analogy is valid also on the technical level: if one forms the dual likelihood ([36]) from score function $Z_t^{\mathrm{MG}}$, the observed information is, indeed, $\langle Z, Z \rangle_t$.

If the drift $Z^{\mathrm{DR}}$ in (4.1) is nonzero, the analysis applies directly only to $Z^{\mathrm{MG}}$. So long as $T$ is small or moderate, however, the variability in $Z^{\mathrm{MG}}$ is the main part of the variability in $Z$. Specifically, both when $T \to 0$ and $T \to +\infty$, $Z_T^{\mathrm{MG}} = O_p(T^{1/2})$ and $Z_T^{\mathrm{DR}} = O_p(T)$. Thus, the bias due to analyzing $Z_t^{\mathrm{MG}}$ in lieu of $Z_t$ becomes a problem only for large $T$. At the same time the present methods provide estimates for the variation in $Z_t^{\mathrm{MG}}$ for small and moderate $T$, whereas the variation in $Z_t^{\mathrm{DR}}$ can only be consistently estimated when $T \to +\infty$ (by Girsanov's theorem). Thus, we recommend our current methods for moderate $T$, while one should use other approaches when dealing with a time span $T$ that is long.

4.5. *Financial applications*: *An instance where variance and variation relate exactly.* It is quite common in finance to encounter the case from Section 4.4, where $Z$ itself is a martingale, or where one is interested in $Z^{\mathrm{MG}}$ only. We here show a conceptual example of this, where one wants to test whether the residual $Z$ is zero, or study the distribution of the residual under the so-called *Risk Neutral* or *Equivalent Martingale Measure* $P^*$.

If $P$ is the true, actual (physical) probability distribution under which data is collected, $P^*$ is, by contrast, a probability measure equivalent to $P$ in the sense of mutual absolute continuity, and it satisfies the condition that the discounted value of all traded securities must be $P^*$-martingales. The values of financial assets, consequently, are expectations under $P^*$. For further details, refer to [8, 9, 11, 24, 25, 28].

If the residual $Z$ relates to the value of a security, one is often interested in its behavior under $P^*$, rather than under $P$. Specifically, we shall see that one is interested in $Z^{\mathrm{MG}*}$, where this is the martingale part in the Doob–Meyer decomposition (4.1), *when taken under* $P^*$:

$$Z_t = Z_0 + Z_t^{\mathrm{DR}*} + Z^{\mathrm{MG}*} \qquad \text{w.r.t. } P^*. \tag{4.2}$$

The quadratic variation $\langle Z, Z \rangle = \langle Z^{\mathrm{MG}}, Z^{\mathrm{MG}} \rangle = \langle Z^{\mathrm{MG}*}, Z^{\mathrm{MG}*} \rangle$ is the same under $P$ and $P^*$, but, under the latter distribution it refers to the behavior of $Z^{\mathrm{MG}*}$ rather than $Z^{\mathrm{MG}}$.

A simple example follows from the motivating application in the Introduction. Suppose that $\Xi$ and $S$ are both discounted securities prices, and that one seeks to offset risk in $\Xi$ by holding $\rho$ units of $S$. The residual is then, itself, the discounted value on the unhedged part of $\Xi$. Under $P^*$, therefore, $Z$ is a martingale, $Z_t = Z_0 + Z_t^{\mathrm{MG}*}$. A deeper example is encountered in [44], where we analyze implied volatilities. In both these cases, in order to put a value on the risk involved in $Z^{\mathrm{MG}*}$, one is interested in bounds on the quadratic variation $\langle Z, Z \rangle$, under $P^*$. This will help, for example, in pricing spread options on $Z$.

How do our results for probability $P$ relate to $P^*$? They simply carry over, unchanged, to this probability distribution. Theorems 1 and 2 remain valid by absolute continuity of $P^*$ under $P$. In the case of limiting results, such as those in Propositions 1 and 3 (in Section 6) and the development for goodness of fit in Section 4, we also invoke the mode of stable convergence (Section 2.8) together with the fact that $dP^*/dP$ is measurable with respect to the underlying $\sigma$-field $\mathcal{F}_T$.

Finally, if one wants to test a null hypothesis $H_0$ that $Z^{\mathrm{MG}*}$ is constant, then $H_0$ is equivalent to asking whether $\langle Z, Z \rangle_T$ is zero (whether under $P$ or $P^*$). This can again be answered with our distributional results above. In the case of the example in this section, the $H_0$ of fully offsetting the risk in $\Xi$ also tests whether $Z$ itself is constant.

**5. Conclusion.** This paper provides a methodology to analyze the association between Itô processes. Under the framework of nonparametric, one-factor regression, we obtain the distributions of estimators of residual variation $\langle Z, Z \rangle$. We then use this in a variety of measures of goodness of fit. We also show how the method yields a procedure to test the appropriateness of a parametric regression model. The limiting distributions identify two sources of uncertainty, one from the discrete nature of the data process, the other from the estimation procedure. Interestingly, among the class of estimators $\widehat{\langle Z, Z \rangle}^{(\alpha)}$ under consideration, $\alpha \in [0, 1]$, discrete-time sampling only impacts the "variance" component. On the other hand, different estimation schemes lead to different biases in the asymptotics.

ANOVA for diffusions permits inference over a time interval. This is because the error terms in the quadratic variation $\widehat{\langle Z, Z \rangle}^{(\alpha)}$ of residuals and, hence, the error terms in the goodness of fit measures, converge as a process, whereas the errors in the estimated regression parameters $\hat{\rho}_t$ are asymptotically independent from one time point to the next. This feature of time aggregation makes ANOVA a natural procedure to determine the adequacy of an adopted model. Also, the ANOVA is better posed in that the rate of the convergence is the square of the rate for $\hat{\rho}_t - \rho_t$.

The "ANOVA for diffusions" approach is appealing also from the position of applications. As long as one can collect the data as a process, one can rely on the proposed ANOVA methodology to draw inference without imposing parametric structure on the underlying process. In financial applications, as in Section 4.5, it can test whether a financial derivative can be fully hedged in another asset. In the event of nonperfect hedging, Theorems 1 and 2 tell us how to quantify the amount of hedging error, as well as its distribution.

**6. Convergence in law—Proofs and further results.** In the following we deal with processes that are exemplified by $[Z, Z] - \langle Z, Z \rangle$. We mostly follow [32].

PROOF OF PROPOSITION 1.    The applicable parts of the proof of the cited Theorems 5.1 and 5.5 of [32] carry over directly under Assumptions A(i) and A(ii). When modifying the proofs, as appropriate, $t_* = \max(t_i^{(n)} \leq t)$ replaces $[tn]/n$, $\delta_n$ replaces $n^{-1}$ and so on. For example, the right-hand side of their equation (5.10) on page 290 becomes $K\delta_n^2$. The main change due to the nonequidistant case occurs in part (iii) of Jacod and Protter's Lemma 5.3, pages 291 and 292, where in the definition of $\alpha_n$, $\frac{t}{2}B_r^{ik}$ should be replaced by $(H(t_r + t) - H(t_r))B_r^{ik}$. Assumptions A(i) and A(ii) are clearly sufficient.    □

Note that the result extends in an obvious fashion to the case of multidimensional $Z = (Z^{(1)}, \ldots, Z^{(p)})$. Also, instead of studying $[Z, Z] - \langle Z, Z \rangle$, one can, like [32], state the (more general) result for $\int_0^t (Z_u^{(i)} - Z_*^{(i)}) \, dZ_u^{(j)}$.

In the sequel we shall also be using a triangular array form of Proposition 1; see the ends of the proofs of Theorems 1 and 2.

PROPOSITION 3 (Triangular array version of the discretization theorem). *Let $Z$ be a vector Itô process for which $\int_0^T \|\langle Z, Z \rangle_t'\|^2 \, dt < \infty$ and $\int_0^T \|\tilde{Z}_t\|^2 \, dt < \infty$ a.s. Also, suppose that $Z^{(n)}$, $i = 1, 2, \ldots$, are Itô processes satisfying the same requirement, uniformly. Suppose that the (vector) Brownian motion $W$ is the same in the Itô process representations of $Z$ and of all the $Z^{(n)}$, that is,*

$$(6.1) \qquad dZ_u^{(n),\mathrm{MG}} = \sigma_u^{(n)} \, dW \quad \text{and} \quad dZ_u^{\mathrm{MG}} = \sigma_u \, dW.$$

*Suppose that*

$$(6.2) \qquad \int_0^T \|\sigma_u^{(n)} - \sigma_u\|^4 \, du = o_p(1).$$

*Then, subject to Assumptions A(i) and (ii), the processes $\frac{1}{\sqrt{\Delta t^{(n)}}} \int_0^t (Z_u^{(i,n)} - Z_*^{(i,n)}) \, dZ_u^{(j,n)}$ converge jointly with the processes $\frac{1}{\sqrt{\Delta t^{(n)}}} \int_0^t (Z_u^{(i)} - Z_*^{(i)}) \, dZ_u^{(j)}$ to the same limit.*

If one requires stable convergence, one just imposes System Assumption III; see Theorem 11.2, page 338, and Theorem 15.2(c), page 496, of [30].

PROOF OF PROPOSITION 3. This is mainly a modification of the development on page 292 and the beginning of page 293 in [32], and the further development in (their) Theorem 5.5 is straightforward. Again we recollect that $H$ [from Assumptions A(i) and (ii)] is Lipschitz continuous.

Note that to match the end of the proof of Theorem 5.1, we really need $o_p(\delta_n^4)$. This, of course, follows by appropriate use of subsequences. □

Finally, we show the result for Section 3.4.

PROOF OF PROPOSITION 2. Let $t_*$ be the largest grid point smaller than or equal to $t$. Then by Itô's lemma, $d(Z_t - Z_{t_*})^4 = 4(Z_t - Z_{t_*})^3 \, dZ_t + 6(Z_t - Z_{t_*})^2 \, d\langle Z, Z \rangle_t$. Hence,

$$
\begin{aligned}
[Z, Z, Z, Z]_t &= \sum_{t_{i+1} \leq t} \left[ 4 \int_{t_i}^{t_{i+1}} (Z_u - Z_{t_i})^3 \, dZ_u \right. \\
(6.3) &\qquad\qquad \left. + 6 \int_{t_i}^{t_{i+1}} (Z_u - Z_{t_i})^2 \, d\langle Z, Z \rangle_u \right] \\
&= 6 \sum_{t_{i+1} \leq t} \int_{t_i}^{t_{i+1}} (Z_u - Z_{t_i})^2 \, d\langle Z, Z \rangle_t + O_p(\delta^{(n)3/2-\varepsilon}),
\end{aligned}
$$

using Lemma 2 below.

Define an interpolated version of $[Z, Z] = [Z, Z]^{\mathcal{G}}$ by letting $[Z, Z]_t^{\text{interpol}} = (Z_t - Z_{t_*})^2 + [Z, Z]_{t_*}$. Then, again by Itô's lemma,

$$
\begin{aligned}
(6.4) \quad & \sum_{t_{i+1} \leq t} \int_{t_i}^{t_{i+1}} (Z_u - Z_{t_i})^2 \, d\langle Z, Z \rangle_t \\
& = \tfrac{1}{4} \langle ([Z, Z]^{\text{interpol}} - \langle Z, Z \rangle), ([Z, Z]^{\text{interpol}} - \langle Z, Z \rangle) \rangle_{t_*}.
\end{aligned}
$$

Putting together (6.3) and (6.4), and Proposition 1, as well as Corollary VI.6.30, page 385 of [33], we obtain (3.8).

Replacing $Z$ by $\widehat{Z}$, and creating interpolated versions of $\widehat{Z}$ and $[\widehat{Z}, \widehat{Z}]$ as at the beginning of the proof of Theorem 1 below, (3.7) also follows. $\quad \square$

## 7. Proofs of main results.

7.1. *Notation.* In the following proofs $t_i$ *always* means $t_i^{(n)}$, $h$ always means $h_n$, $\overline{\Delta t}$ means $\overline{\Delta t}^{(n)}$ and so on. Also, $t_*$ is the largest grid point less than $t$, that is, $t_* = \max\{t_i^{(n)} \leq t\}$. We sometimes write $\langle X, X \rangle_t$ as $\langle X \rangle_t$, and $\langle X, X \rangle_t'$ as $\langle X \rangle_t'$ for simplicity. Also, for convenience, we adopt the following shorthand for smoothness assumptions for Itô processes:

ASSUMPTION B (*Smoothness*).
B.0($X$): $X$ is in $C^1[0, T]$.
B.1($X, Y$): $\langle X, Y \rangle_t$ is in $C^1[0, T]$.
B.2($X$): the drift part of $X$ ($X^{\text{DR}}$) is in $C^1[0, T]$.

Assumption B($X$) is equivalent to B.1($X, X$) and B.2($X$).
We shall also be using the following notation:

$$
\Upsilon^X(h) = \sup_{u,s} |X_u - X_s| \quad \text{and} \quad \Upsilon^{XY}(h) = \sup_{u,s} |\langle X, Y \rangle_u' - \langle X, Y \rangle_s'|,
$$

where $\sup_{u,s}$ means $\sup_{u,s \in [0,T]:|u-s| \leq h}$.

7.2. *Preliminary lemmas and proof of Theorem* 1.

LEMMA 1. *Let* $M_{i,n}(t)$, $0 \leq t \leq T$, $i = 1, \ldots, k_n$, $k_n = O(\overline{\Delta t}^{-1})$, *be a collection of continuous local martingales. Suppose that* $\sup_{1 \leq i \leq k_n} \langle M_{i,n}, M_{i,n} \rangle_T = O_p(\overline{\Delta t}^{\beta})$. *Then, for any* $\varepsilon > 0$, $\sup_{1 \leq i \leq k_n} \sup_{0 \leq t \leq T} |M_{i,n}(t)| = O_p(\overline{\Delta t}^{\beta/2 - \varepsilon})$.

The above follows from Lenglart's inequality. As a corollary, the following is true.

LEMMA 2.  *Let $X$ and $Y$ be Itô processes. If $X$ satisfies that $|\widetilde{X}|$ and $\langle X, X \rangle'$ are bounded a.s., then for any $\varepsilon > 0$, $\Upsilon^X(\eta) = O_p((\eta + \overline{\Delta t})^{1/2 - \varepsilon})$. Similarly, if $X$ and $Y$ satisfy that $\widetilde{D}^{XY}$ and $\langle R^{XY}, R^{XY} \rangle$ are bounded (a.s.), then for any $\varepsilon > 0$, $\Upsilon^{XY}(\eta) = O_p((\eta + \overline{\Delta t})^{1/2 - \varepsilon})$.*

LEMMA 3.  *Suppose $X$, $Y$ and $Z$ are Itô processes, and make Assumptions A(i), B.0($Y$), B[($X$), ($Z$)] and B.1($X, Z$). Also, assume that for any $\varepsilon > 0$, as $n \to \infty$, $(\overline{\Delta t})^{1/2 - \varepsilon} h^{-1/2} = o(1)$. Then*

$$\sup_t \frac{1}{h^2} \left| \sum_{t - h \leq t_i < t_{i+1} \leq t} \left[ \int_{t_i}^{t_{i+1}} (X_u - X_{t_i})(Z_u - Z_{t_i}) \, dY_u - \frac{(\Delta t_i)^2}{2} \langle X, Z \rangle'_{t_i} \widetilde{Y}_{t_i} \right] \right|$$

$$= o_p \left( \frac{\overline{\Delta t}}{h} \right)$$

*as $n \to \infty$, where $\widetilde{Y}_u = dY_u / du$.*

PROOF OF LEMMA 3.  Without loss of generality, it is enough to show the result for $X = Z$. This is because one can prove the results for $X$, $Z$ and $X + Z$, and then proceed via the polarization identity. The conditions imposed also mean that the assumptions of Lemma 3 are also satisfied for $X + Z$. Let

$$I_t = \frac{1}{h^2} \sum_{t - h \leq t_i < t_{i+1} \leq t} \left[ \int_{t_i}^{t_{i+1}} (\langle X \rangle_u - \langle X \rangle_{t_i}) \, dY_u - \frac{1}{2} \langle X \rangle'_{t_i} \widetilde{Y}_{t_i} (\Delta t_i)^2 \right],$$

$$II_t = \frac{1}{h^2} \sum_{t - h \leq t_i < t_{i+1} \leq t} \int_{t_i}^{t_{i+1}} [(X_u - X_{t_i})^2 - (\langle X \rangle_u - \langle X \rangle_{t_i})] \, dY_u.$$

Obviously $\sup_t |I_t| = o_p(\overline{\Delta t} h^{-1})$. It is then enough to show

$$\sup_t |II_t| = o_p(\overline{\Delta t} h^{-1}),$$

as follows.

By Itô's lemma,

$$II_t = \underbrace{\frac{2}{h^2} \sum_{t - h \leq t_i < t_{i+1} \leq t} \int_{t_i}^{t_{i+1}} \int_{t_i}^u (X_v - X_{t_i}) \, dX_v^{\mathrm{DR}} \, dY_u}_{II_{t,1}}$$

$$+ \underbrace{\frac{2}{h^2} \sum_{t - h \leq t_i < t_{i+1} \leq t} \int_{t_i}^{t_{i+1}} \left[ \int_{t_i}^u (X_v - X_{t_i}) \, dX_v^{\mathrm{MG}} \right] dY_u}_{II_{t,2}}.$$

Recall that $dX_v^{\mathrm{DR}} = \widetilde{X}_v\, dv$. Then by B.0($Y$), B.2($X$) and the continuity of $X$, one gets

$$|II_{t,1}| \leq \sup_{0 \leq u \leq t} |\widetilde{X}_u| \sup_{0 \leq u \leq t} |\widetilde{Y}_u| \frac{2}{h^2} \sum_{t-h \leq t_i < t_{i+1} \leq t} \int_{t_i}^{t_{i+1}} \left( \int_{t_i}^{u} |X_v - X_{t_i}|\, dv \right) du$$

$$= o_p\!\left( \frac{\delta^{(n)}}{h} \right).$$

For $II_{t,2}$, using integration by parts,

$$\frac{1}{h^2} \sum_{t-h \leq t_i < t_{i+1} \leq t} \int_{t_i}^{t_{i+1}} \left[ \int_{t_i}^{u} (X_v - X_{t_i})\, dX_v^{\mathrm{MG}} \right] du$$

$$(7.1) \qquad = \frac{1}{h^2} \sum_{t-h \leq t_i < t_{i+1} \leq t} \int_{t_i}^{t_{i+1}} (X_u - X_{t_i})[(\Delta t_i) - (u - t_i)]\, dX_u^{\mathrm{MG}}.$$

Equation (7.1) has quadratic variation bounded by

$$\frac{1}{h^4} \sup_i \sum_{t-h \leq t_i < t_{i+1} \leq t} \int_{t_i}^{t_{i+1}} (X_u - X_{t_i})^2 (t_{i+1} - u)^2\, d\langle X \rangle_u$$

$$(7.2) \qquad \leq \frac{1}{h^4} \sup_i \langle X \rangle_{t_i}' \sup_t \sum_{t-h \leq t_i < t_{i+1} \leq t} \int_{t_i}^{t_{i+1}} (\Upsilon^X(\delta^{(n)}))^2 (t_{i+1} - u)^2\, du$$

$$= O_p(\overline{\Delta t}^{3-\varepsilon} h^{-3})$$

by Lemma 2 under Assumptions A(i) and B($X$). Following Lemma 1 and B.0($Y$), $\sup_t |II_{t,2}| = O_p(\overline{\Delta t}^{3/2-\varepsilon} h^{-3/2})$. The result then follows.

DEFINITION 5.  Suppose $X$ and $Y$ are continuous Itô processes. Let

$$(7.3) \qquad B_{1,i,t}^{XY} = \begin{cases} \frac{1}{h} \int_{t_i-h}^{t \wedge t_i} ((t_i - h) - u)\, d\langle X, Y \rangle_u', & t \geq t_i - h, \\ 0, & \text{otherwise} \end{cases}$$

and

$$(7.4)\ B_{2,i,t}^{XY} = \begin{cases} \frac{[2]}{h} \sum_{t_i-h \leq t_j \leq t_{j+1} \leq t_i \wedge t} \int_{t_j}^{t_{j+1}} (X_s - X_{t_j})\, dY_s, & t \geq t_i - h, \\ 0, & \text{otherwise,} \end{cases}$$

where [2] indicates symmetric representation s.t. $[2] \int X\, dY = \int X\, dY + \int Y\, dX$.

Note that by integration by parts via Itô's lemma, $B_{1,i,t_i}^{XY} = \frac{1}{h}(\langle X, Y \rangle_{t_i} - \langle X, Y \rangle_{t_i-h}) - \langle X, Y \rangle_{t_i}'$ and, hence, $\widehat{\langle X, Y \rangle}_{t_i}' - \langle X, Y \rangle_{t_i}' = B_{1,i,t_i}^{XY} + B_{2,i,t_i}^{XY}$.

LEMMA 4. *Under Assumptions* A(i) *and* B$(X, Y, \langle X, Y \rangle')$ *and the order selection of* $h^2 = O(\overline{\Delta t})$, *for any* $\varepsilon > 0$,

$$\sup_{0 \le t_i \le T} |B_{1,i,t_i}^{XY}| = O_p(\overline{\Delta t}^{1/4-\varepsilon}) \quad and \quad \sup_{0 \le t_i \le T} |B_{2,i,t_i}^{XY}| = O_p(\overline{\Delta t}^{1/4-\varepsilon}).$$

*In particular,*

$$\sup_{0 \le t_i \le T} |\widehat{\langle X, Y \rangle}_{t_i}' - \langle X, Y \rangle_{t_i}'| = O_p(\overline{\Delta t}^{1/4-\varepsilon}).$$

*In addition, under condition* $B(\langle R^{XY}, R^{ZV} \rangle')$,

$$(7.5) \qquad \sup_{t_i} \left| \langle B_{1,i}^{XY}, B_{1,i}^{ZV} \rangle_{t_i} - \frac{h}{3} \langle R^{XY}, R^{ZV} \rangle'_{t_i} \right| = O_p(h^{3/2-\varepsilon}),$$

$$(7.6) \qquad \sup_{t} |\langle B_{2,i}^{XY}, B_{2,i}^{ZV} \rangle_t - G_t| = o_p(\overline{\Delta t}^{1/2}),$$

*where* $G_t = \frac{1}{h^2} \sum_{t-h \le t_i < t_{i+1} \le t} (\langle X, Z \rangle'_{t_i} \langle Y, V \rangle'_{t_i} + \langle X, V \rangle'_{t_i} \langle Y, Z \rangle'_{t_i})(\Delta t_i)^2$ *and*

$$(7.7) \qquad \sup_{t_i} |\langle B_{1,i}^{ZV}, B_{2,i}^{XY} \rangle_{t_i}| = O_p\left(\frac{\overline{\Delta t}}{\sqrt{h}}\right).$$

PROOF OF LEMMA 4. Using Definition 2, write

$$(7.8) \quad B_{1,i,t}^{XY} = \underbrace{\frac{1}{h} \int_{t_i-h}^{t} ((t_i - h) - u) \, dR_u^{XY}}_{B_{1,i,t}^{XY,\mathrm{MG}}} + \underbrace{\frac{1}{h} \int_{t_i-h}^{t} ((t_i - h) - u) \, dD_u^{XY}}_{B_{1,i,t}^{XY,\mathrm{DR}}}.$$

Under Assumption B.2($\langle X, Y \rangle'$), $\sup_i |B_{1,i,t_i}^{XY,\mathrm{DR}}| = O_p(h)$. Also, by Assumptions A(i) and B.1($R^{XY}, R^{XY}$), we have that $\sup_{0 \le u \le T} \langle R^{XY} \rangle'_u = O_p(1)$, whence $\sup_i \langle B_{1,i}^{XY}, B_{1,i}^{XY} \rangle_T = O_p((\overline{\Delta t}^{(n)})^{1/2})$. So $\sup_i \sup_t |B_{1,i,t}^{XY,\mathrm{MG}}| = O_p(\overline{\Delta t}^{1/4-\varepsilon})$ by Lemma 1 and, hence, $\sup_i |B_{1,t_i}^{XY}| = O_p(\overline{\Delta t}^{1/4-\varepsilon})$. By similar methods, Lemmas 1 and 2 can be used to show that $\sup_i |B_{2,t_i}^{XY}| = O_p(\overline{\Delta t}^{1/4-\varepsilon})$.

The orders (7.5)–(7.7) follow from the representation (7.8), and derivations similar to those of the proof of Lemma 3.  □

The following is now immediate from Lemma 4, by Taylor expansion.

COROLLARY 1. *Suppose* $\Xi$, $S$ *and* $\rho$ *are Itô processes, where* $\rho$ *and* $\hat{\rho}$ *are as defined in Section 2. Then, under conditions* A(i), B$(\Xi, S, \langle \Xi, S \rangle', \langle S, S \rangle')$ *and* (2.7), *for any* $\varepsilon > 0$, *we have*

$$\sup_{t_i \in [0,T]} |\hat{\rho}_{t_i} - \rho_{t_i}| = O_p(\overline{\Delta t}^{1/4-\varepsilon}).$$

In the rest of this subsection we shall set, by convention, $\hat{\rho}_t$ to the value $\hat{\rho}_{t*}$, even if the definition in Section 2.4 permits other values of $\hat{\rho}_t$ between sampling points. This is no contradiction since we only use the $\hat{\rho}_{t_i}$ in our definition of $\hat{Z}$ and in the rest of our development.

Similarly, we extend the definition of $\hat{Z}_t$ given at the beginning of Section 2.5 to the case where $t$ is not a sampling point. Specifically,

$$\hat{Z}_t = \hat{Z}_{t_*} + \Xi_t - \Xi_{t*} - \hat{\rho}_{t*}(S_t - S_{t*})$$

$$= \hat{Z}_{t_*} + \Xi_t - \Xi_{t*} - \int_{t*}^t \hat{\rho}_u \, dS_u$$

by our convention for $\hat{\rho}_t$. We emphasize that $\hat{Z}_t$ is no more observable than $S_t$ or $\Xi_t$. By this definition of $\hat{Z}$, in view of (1.6) and (1.5),

$$(7.9) \qquad \langle \hat{Z}, \hat{Z} \rangle_t = \langle Z, Z \rangle_t + \int_0^t (\hat{\rho}_u - \rho_u)^2 \, d\langle S, S \rangle_u.$$

We then obtain the following preliminary result for Theorem 1.

PROPOSITION 4. *Under the conditions of Theorem* 1, $(\langle \hat{Z}, \hat{Z} \rangle_t - \langle Z, Z \rangle_t)/$ $\sqrt{\overline{\Delta t}} = D_t + o_p(1)$, *uniformly in* $t$, *where* $D_t$ *is given by* (3.2).

PROOF. (i) Let $L_{i,n}(t) = \sum_{j=1}^2 L_{j,i,n}(t)$, where for $j = 1, 2$, $L_{j,i,n}(t) = (B_{j,t \wedge t_i}^{\Xi S} - \rho_{t_i-h} B_{j,t \wedge t_i}^{SS})/\langle S, S \rangle_{t_i-h}'$ for $t \geq t_i - h$, and zero otherwise. It is now easily seen from Lemmas 2 and 4, and Corollary 1, and by the same methods as in these results, that $L_{i,n}(t)$ approximates $\hat{\rho}_t - \rho_t$ and, in particular, that

$$(7.10) \quad \int_0^t (\hat{\rho}_u - \rho_u)^2 \, d\langle S, S \rangle_u = \sum_{t_{i+1} \leq t} L_{i,n}^2(t_i) \langle S, S \rangle_{t_i-h}' \Delta t_i + O_p(\overline{\Delta t}^{3/4-\varepsilon})$$

uniformly in $t$, for any $\varepsilon > 0$. We now show that

$$(7.11) \begin{aligned} \sum_{t_{i+1} \leq t} L_{i,n}^2(t_i) \langle S, S \rangle_{t_i-h}' \Delta t_i \\ = \sum_{t_{i+1} \leq t} \langle L_{i,n} \rangle_{t_i} \langle S, S \rangle_{t_i-h}' \Delta t_i + O_p(\overline{\Delta t}^{3/4-\varepsilon}). \end{aligned}$$

To this end, set $Y^{(i)}(t) = L_{i,n}^2(t) - \langle L_{i,n} \rangle_t$ and

$$(7.12) \quad Z_{n,t} = \sum_{t_{i+1} \leq t} Y^{(i)}(t_i) \langle S, S \rangle_{t_i-h}' \Delta t_i + Y^{(i)}(t_*) \langle S, S \rangle_{t_*-h}' \Delta t_*.$$

By Lenglart's inequality ([33], Lemma I.3.30, page 35), (7.11) follows if one can show that $\langle Z_n, Z_n \rangle_T = O_p(\overline{\Delta t}^{3/2-\varepsilon})$, for any $\varepsilon > 0$. This is what we do in the rest of (i).

Since $\langle L_{i,n}, L_{j,n} \rangle_t = 0$ if $(t_i - h, t_i]$ and $(t_j - h, t_j]$ are disjoint for any $t_i \leq T$ and $t_j \leq T$, the quadratic variation of the $Z_n$ is considered in the overlapping time interval:

$$\langle Z_n, Z_n \rangle_T$$

$$\leq \sup_{0 \leq t \leq T} (\langle S, S \rangle_t')^2 \left| \sum_i \sum_j (\Delta t_i)(\Delta t_j) \int_{(t_i - h, t_i] \cap (t_j - h, t_j]} d\langle Y^{(i)}, Y^{(j)} \rangle_u \right|$$

$$(7.13) \quad \leq 2 \sup_{0 \leq t \leq T} (\langle S, S \rangle_t')^2 \sum_{i \leq j} (\Delta t_i)(\Delta t_j) I_{\{i < j : t_i > t_j - h\}} \left( \int_{(t_j - h, t_i]} d\langle Y^{(i)} \rangle_u \right)^{1/2}$$

$$\times \left( \int_{(t_j - h, t_i]} d\langle Y^{(j)} \rangle_u \right)^{1/2}.$$

The above line follows from the Kunita–Watanabe inequality (page 69 of [38]).

By Itô's formula, for $t_i - h < t < t_i$, $\langle Y^{(i)} \rangle_t' = 4 L_{i,n}^2(t) \langle L_{i,n} \rangle_t'$. Hence, by Lemma 4, $\langle Y^{(i)} \rangle_t' \leq U_1 \langle L_{i,n} \rangle_t'$, where $U_1$ is defined independently of $i$, and $U_1 = O_p(\overline{\Delta t}^{1/2 - \varepsilon})$. Also, by the Kunita–Watanabe inequality and the Cauchy–Schwarz inequality,

$$(7.14) \quad \langle L_{i,n} \rangle_t' \leq 4 \sup_{0 \leq u \leq T} \frac{1}{(\langle S, S \rangle_u')^2} \sum_{j=1}^{2} (\langle B_{j,i}^{\Xi S} \rangle_t' + \rho_{t_i - h}^2 \langle B_{j,i}^{SS} \rangle_t'),$$

where $B_{j,i,t}^{XY}$, $j = 1, 2$, is defined before Lemma 4. Thus, it is enough that (7.13) is $O_p(\overline{\Delta t}^{3/2 - \varepsilon})$ in the two cases where $\langle Y^{(i)} \rangle_t'$ is replaced by (a) $U_1 \langle B_{1,i}^{XY} \rangle_t'$, for $(X, Y) = (\Xi, S)$ and $(S, S)$, and (b) $U_1 \langle B_{2,i}^{XY} \rangle_t'$, for $(X, Y) = (\Xi, S), (S, \Xi)$ and $(S, S)$. We do this in the case of $\langle B_{1,i}^{XY} \rangle_t'$. The second case is similar.

Obviously, on $t_i - h < t < t_i$, $\langle B_{1,i}^{XY} \rangle_t' = \frac{1}{h^2}(t - (t_i - h))^2 \langle R^{XY} \rangle_t'$. Also, set $N_n = \sup_t \#\{j : t \leq t_j \leq t + h\}$ and $\delta_-^{(n)} = \min (t_{i+1}^{(n)} - t_i^{(n)})$, and note that $N_n = O(h / \delta_-^{(n)}) = O(h / \overline{\Delta t})$ under Assumption A. Since $\sup_{0 \leq u \leq T} \langle R^{XY} \rangle_u' < \infty$, the part of equation (7.13) which is attributable to the $B_1^{XY}$ term becomes, up to an $O_p(1)$ factor,

$$(7.15) \quad U_1 \frac{\overline{\Delta t}^2}{h^2} \sum_{i \leq j : t_j - t_i < h} [h^3 - (t_j - t_i)^3]^{1/2} [t_i - (t_j - h)]^{3/2} = O_p(\overline{\Delta t}^{3/2 - \varepsilon}).$$

(ii) To finish the proof of the proposition, we wish to show that, as $(\overline{\Delta t}/h)^{1/2} \to c$,

$$(7.16) \quad \overline{\Delta t}^{-1/2} \sum_{t_{i+1} \leq t} \langle L_{i,n} \rangle_{t_i} \langle S, S \rangle_{t_i - h}' \Delta t_i \to D_t$$

in probability, for each $t$. (This is enough, since the convergence of increasing functions to an increasing function is automatically uniform.) Since $\sup_i |\langle L_{1,i,n}, L_{2,i,n} \rangle_{t_i}| = O_p(\overline{\Delta t}^{3/4})$ from Lemma 4, it follows that it is sufficient to prove separately that

$$
\begin{aligned}
(7.17) \qquad & \overline{\Delta t}^{-1/2} \sum_{t_{i+1} \leq t} \langle L_{1,i,n} \rangle_{t_i} \langle S, S \rangle'_{t_i - h} \Delta t_i \\
& \xrightarrow{P} \frac{1}{3c} \int_0^t \langle \rho, \rho \rangle'_u \, d\langle S, S \rangle_u
\end{aligned}
$$

and

$$
\begin{aligned}
(7.18) \qquad & \overline{\Delta t}^{-1/2} \sum_{t_{i+1} \leq t} \langle L_{2,i,n} \rangle_{t_i} \langle S, S \rangle'_{t_i - h} \Delta t_i \\
& \xrightarrow{P} c \int_0^t \left( \frac{\langle \Xi, \Xi \rangle'_u}{\langle S, S \rangle'_u} - \rho_u^2 \right) H'(u) \, d\langle S, S \rangle_u.
\end{aligned}
$$

Equation (7.17) follows directly from the approximation (7.5) in Lemma 4. It remains to show (7.18). This is what we do for the rest of the proof.

Let $A$ be an Itô process, which we shall variously take to be $\frac{1}{\langle S, S \rangle'}$, $-\frac{2\rho}{\langle S, S \rangle'}$ and $\frac{\rho^2}{\langle S, S \rangle'}$. Consider a subproblem of (7.18), that of the convergence of

$$
(7.19) \qquad \overline{\Delta t}^{-1/2} \sum_{t_{i+1} \leq t} \langle B_{2,i}^{XY}, B_{2,i}^{ZV} \rangle_{t_i} A_{t_i - h} \Delta t_i.
$$

By (7.6) in Lemma 4, this is (uniformly in $t$) equal to

$$
\overline{\Delta t}^{-1/2} h^{-2} \sum_{t_{i+1} \leq t} \sum_{t_i - h \leq t_j \leq t_{j+1} \leq t_i} f(t_j) (\Delta t_j)^2 A_{t_i - h} \Delta t_i + o_p(1),
$$

where $f(t) = \langle X, Z \rangle'_t \langle Y, V \rangle'_t + \langle X, V \rangle'_t \langle Y, Z \rangle'_t$. By interchanging the two summations, (7.19) then becomes, up to $o_p(1)$,

$$
\begin{aligned}
& \overline{\Delta t}^{-1/2} h^{-1} \sum_{t_{j+2} \leq t} f(t_j) (\Delta t_j)^2 h^{-1} \sum_{t_{j+1} \leq t_i \leq t_j + h} A_{t_i - h} \Delta t_i \\
& = \overline{\Delta t}^{-1/2} h^{-1} \sum_{t_{j+2} \leq t} f(t_j) A_{t_j} (\Delta t_j)^2.
\end{aligned}
$$

This results because the difference between the last two terms is bounded by

$$
\begin{aligned}
\overline{\Delta t}^{-1/2} h^{-1} \sum_{t_{j+2} \leq t} |f(t_j)| (\Delta t_j)^2 \Upsilon^A(h) & \leq \sup_t |f(t)| \Upsilon^A(h) \overline{\Delta t}^{1/2} h^{-1} H_n(t) \\
& = o_p(1),
\end{aligned}
$$

by Lemma 2. Hence, (7.19) converges to $c \int_0^t f(u) A_u \, dH(u) = c \int_0^t f(u) A_u \times H'(u) \, du$ by Assumption $A$ and since $A$ and $f$ are bounded and continuous. Note that $H$ is absolutely continuous since it is Lipschitz. The result (7.18) now follows by aggregating this convergence for the cases of $\langle B_{2,i}^{\Xi S}, B_{2,i}^{\Xi S} \rangle$ $(A = \frac{1}{\langle S,S \rangle'})$, $\langle B_{2,i}^{\Xi S}, B_{2,i}^{SS} \rangle$ $(A = -\frac{2\rho}{\langle S,S \rangle'})$ and $\langle B_{2,i}^{SS}, B_{2,i}^{SS} \rangle$ $(A = \frac{\rho^2}{\langle S,S \rangle'})$.   □

PROOF OF THEOREM 1.  In view of Proposition 4, it is enough to show that

(7.20)    $$\sup_t |[\widehat{Z}, \widehat{Z}]_t - \langle \widehat{Z}, \widehat{Z} \rangle_t - ([Z, Z]_t - \langle Z, Z \rangle_t)| = o_p(\overline{\Delta t}^{-1/2}).$$

From (7.9), $\langle \widehat{Z}, \widehat{Z} \rangle_t' - \langle Z, Z \rangle_t' = (\hat{\rho}_t - \rho_t)^2 \langle S, S \rangle'$, and similarly for the drift of $\widehat{Z}$ and $Z$. The result then follows from Proposition 3 and Corollary 1.   □

**7.3.** *Additional lemmas for Theorem 2, and proof of the theorem.*  By the same methods as above, we obtain the following two lemmas, where the first is the key step in the second.

LEMMA 5.  *Let $X$, $Y$ and $A$ be Itô processes. Let $h = O(\overline{\Delta t}^{1/2})$. Assume Assumptions A, B.1[$(X, X), (A, A), (Y, Y), (A, X), (A, Y)$] and B.2[$(X), (A), (Y)$]. Define*

$$V_t^{XY} = \sum_{t_{i+1} \leq t} (\Delta X_{t_i})(\Delta Y_{t_i}) + (X_t - X_{t_*})(Y_t - Y_{t_*}),$$

(7.21)    $$U(t) = \frac{1}{h} \int_0^t \sum_i A_u I_{(t_i, t_{i+1} \in [u-h, u])} \Delta V_{t_i} \, du$$

$$- \frac{1}{h} \sum_{t_{i+1} \leq t} A_{t_i} \Delta V_{t_i} (h - \Delta t_i).$$

*Then $\sup_t |U(t)| = o_p(\overline{\Delta t}^{1/2})$.*

LEMMA 6.  *Let $\Xi$, $S$, $\rho$ and $Z$ be the Itô processes defined in earlier sections. Subject to the regularity conditions in Lemma 5 with $(X, Y) = (\Xi, S)$, $(S, S)$, or $(\rho, S)$, and with $A = \rho$, $\rho^2$ or $Z$,*

$$\int_0^t (\hat{\rho}_u - \rho_u) \, d\langle \Xi, S \rangle_u = [\Xi, Z]_t - [Z, Z]_t$$

$$- \frac{h}{3} \int_0^t \rho_u \, d\langle \rho, \langle S, S \rangle' \rangle_u + o_p(h),$$

*uniformly in $t$.*

PROOF OF THEOREM 2. By definition, $\Delta \langle \widetilde{Z, Z} \rangle_{t_i} = \Delta [\Xi, \widehat{Z}]_{t_i}$. Since $d\langle \Xi, Z \rangle_t = d\langle Z, Z \rangle_t$ by assumption, and by subtracting and adding $\langle \Xi, \widehat{Z} \rangle_t$,

$$\overline{\Delta t}^{-1/2}(\langle \widetilde{Z, Z} \rangle_t - \langle Z, Z \rangle_t)$$
$$= \underbrace{\overline{\Delta t}^{-1/2}([\Xi, \widehat{Z}]_t - \langle \Xi, \widehat{Z} \rangle_t)}_{C_1} + \underbrace{\overline{\Delta t}^{-1/2}(\langle \Xi, \widehat{Z} \rangle_t - \langle \Xi, Z \rangle_t)}_{C_2}.$$

First notice that

$$(7.22) \qquad \langle \Xi, \widehat{Z} \rangle'_t = \langle \Xi, Z \rangle'_t + (\rho_t - \hat{\rho}_t)\langle \Xi, S \rangle'_t.$$

Also, as $\overline{\Delta t}^{1/2}/h \to c$, (7.22) and Lemma 6 show that

$$C_2 = \frac{1}{\sqrt{\overline{\Delta t}^{(n)}}} \int_0^t (\rho_u - \hat{\rho}_u) \, d\langle \Xi, S \rangle_u$$
$$= \frac{[Z, Z]_t - [\Xi, Z]_t}{\sqrt{\overline{\Delta t}^{(n)}}} + \frac{1}{3c} \int_0^t \rho_u \, d\langle \rho, \langle S, S \rangle' \rangle_u + o_p(1) \qquad \text{uniformly in } t.$$

Since $\langle \Xi, Z \rangle_t = \langle Z, Z \rangle_t$, it follows that

$$\overline{\Delta t}^{1/2}(\langle \widetilde{Z, Z} \rangle_t - \langle Z, Z \rangle_t)$$
$$= \overline{\Delta t}^{1/2}([\Xi, \widehat{Z}]_t - \langle \Xi, \widehat{Z} \rangle_t + [Z, Z]_t - [\Xi, Z]_t)$$
$$\quad + \frac{1}{3c} \int_0^t \rho_u \, d\langle \rho, \langle S, S \rangle' \rangle_u + o_p(1)$$
$$= \overline{\Delta t}^{1/2}([Z, Z]_t - \langle Z, Z \rangle_t) + \frac{1}{3c} \int_0^t \rho_u \, d\langle \rho, \langle S, S \rangle' \rangle_u$$
$$\quad + \overline{\Delta t}^{1/2}([\Xi, \widehat{Z} - Z]_t - \langle \Xi, \widehat{Z} - Z \rangle_t) + o_p(1).$$

The last component in the above, $\overline{\Delta t}^{1/2}([\Xi, \widehat{Z} - Z]_t - \langle \Xi, \widehat{Z} - Z \rangle_t)$, goes to zero in probability by Proposition 3, since $\langle \Xi, \Xi \rangle'_u$, $\langle \widehat{Z} - Z, \widehat{Z} - Z \rangle'_u$, and $\langle \Xi, \widehat{Z} - Z \rangle'_u$ satisfy the conditions of this proposition. This is in view of Corollary 1. The argument is similar to that at the end of the proof of Theorem 1.  □

**Acknowledgments.** We would like to thank the referees for valuable suggestions. The authors would like to thank The Bendheim Center for Finance at Princeton University, where part of the research was carried out.

Department of Statistics
University of Chicago
5734 University Avenue
Chicago, Illinois 60637
USA
E-mail: mykland@galton.uchicago.edu

Department of Statistics
Carnegie Mellon University
Pittsburgh, Pennsylvania 15213
USA
E-mail: lzhang@stat.cmu.edu
and
Department of Finance
University of Illinois at Chicago
Chicago, Illinois 60607
USA
E-mail: lanzhang@uic.edu